\documentclass[reprint,
 amsmath,amssymb,
aps,
pre,
]{revtex4-2}

\usepackage{graphicx,subfigure}
\usepackage{amsfonts,amsbsy,amscd,amsgen,amsthm,dsfont,amsmath,amssymb,mathrsfs,mathtools,array}
\usepackage[]{natbib}
\usepackage[colorlinks=true,allcolors=blue]{hyperref}
\usepackage{bm,multirow,esvect,cancel}
\theoremstyle{definition}
\usepackage{tabularx}
\usepackage{xcolor,comment,soul}
\usepackage{enumitem}
\usepackage[utf8]{inputenc}
\usepackage[T1]{fontenc}
\usepackage[]{algorithm}
\usepackage{algpseudocode} %automatically loads algorithmicx
\newcommand{\id}{\mathrm d}
\newcommand{\vc}{\bm}

\newcommand{\bxi}{\pmb{\xi}}

\renewcommand{\tilde}{\widetilde}

\DeclareMathAlphabet\mathbfcal{OMS}{cmsy}{b}{n}

\newcommand{\R}{\mathbb{R}}

\renewcommand{\d}{\partial}

 \newcommand{\s}{\sigma}

\newtheorem{thm}{Theorem}
\newtheorem{rem}{Remark}
\newtheorem{lem}{Lemma}
\newtheorem{defn}{Definition}

\newtheorem{cor}{Corollary}

\begin{document}
\title{Mitigation of rare events in multistable systems driven by correlated noise}
\author{Konstantinos Mamis} 
\affiliation{Mathematical Modeling and Applications Laboratory, Hellenic Naval Academy, Chatzikiriakou Avenue, Piraeus, Attica 18539, Greece}
\author{Mohammad Farazmand}\thanks{Corresponding author's email address: farazmand@ncsu.edu}
\affiliation{Department of Mathematics, North Carolina State University,
	2311 Stinson Drive, Raleigh, NC 27695-8205, USA}
\date{\today}

\begin{abstract}
	We consider rare transitions induced by colored noise excitation in multistable systems. 
	We show that undesirable transitions  can be mitigated by 
	a simple time-delay feedback control if the control parameters are judiciously chosen.
	We devise a parsimonious method for selecting the optimal control parameters, without requiring any Monte Carlo simulations of the system. This method relies on a new nonlinear Fokker--Planck equation whose stationary response distribution is approximated by a rapidly convergent iterative algorithm. In addition, our framework allows us to accurately predict, and subsequently suppress, the modal drift and tail inflation in the controlled stationary distribution. We demonstrate the efficacy of our method on two examples, including an optical laser model perturbed by multiplicative colored noise.
\end{abstract}

\maketitle

\section{Introduction} \label{sec:Introduction}
Noise-induced transitions are observed in many areas of science and engineering, such as climatology~\cite{Dakos2008,Lenton2008,Mendez2020},
laser technology \cite{Zhu2010}, ecosystems \cite{Ridolfi2011, Spanio2017, Zeng2017b}, oncology \cite{Bose2011, Idris2016, Yang2014, Zeng2010}, neural systems \cite{Li2016, Li2018}, material science~\cite{Bose2012, Chattopadhyay2016}, turbulence~\cite{Gayout2021, Dallas2020, VanKan2019, Shukla2016,faraz_adjoint,Farazmand2017, farazmand2019b} and thermoacoustics~\cite{Zhang2020}. Although these transitions are often rare, their occurrence may have devastating consequences~\cite{farazmand2019a}. Here, we investigate the ability of time-delay feedback control to mitigate such undesirable transitions.

We focus particularly on stochastically excited multistable dynamical systems. In absence of noise, the equilibria of these systems 
are stable fixed points.
Stochastic excitations, however, instigate rare transitions between these equilibria. We assume that one of the equilibria is desirable
and design a control strategy that mitigates transitions away from it.
Following Farazmand~\cite{Farazmand2020}, we consider a class of time-delay feedback controllers. The time delay, although small, is nonzero
in order to model the delay that occurs in applications between observing the system and actuating the controller.

The crucial difference between the present study and Ref.~\cite{Farazmand2020} is the nature of the noise.
For mathematical convenience, stochastic excitations are usually modeled by delta-correlated white noise.
However, environmental noises in reality have a finite correlation time and therefore cannot be modeled as white noise (see, e.g., Refs.~\cite[Sec. 8.1]{Horsthemke2006} and \cite[Sec. 5.4.1]{Pugachev2001}). 
Here, we consider this more realistic case where the noise is colored and therefore has a finite correlation time.

As shown in Ref.~\cite{Farazmand2020}, for the white noise excitation, the time-delay feedback control has two competing effects. 
One is the deepening of the effective potential well around the desirable equilibrium and hence hindering transitions away from it.
At the same time, the control intensifies the effective noise, facilitating large stochastic excursions. We show that the same competing factors are also operative in the case of colored noise.
As such, choosing the suitable control parameters is a delicate balancing act.

Our main goal is to determine the optimal control parameters that minimize the probability of transitions away from
the desirable equilibrium. In other words, the stationary probability density function (PDF) of the controlled system should be unimodal and concentrated 
around the desirable equilibrium. In principle, the optimal control parameters can be determined by direct Monte Carlo simulations of
the controlled system. However, these simulations are computationally expensive and therefore impractical.

Here, we propose a series of approximations that facilitates a parsimonious estimation of the 
stationary response PDF of the controlled system. This in turn allows us to sweep the control parameter space 
and determine the optimal control parameters in a computationally inexpensive manner.

\begin{figure*}
	\centering
	\includegraphics[width=0.75\textwidth]{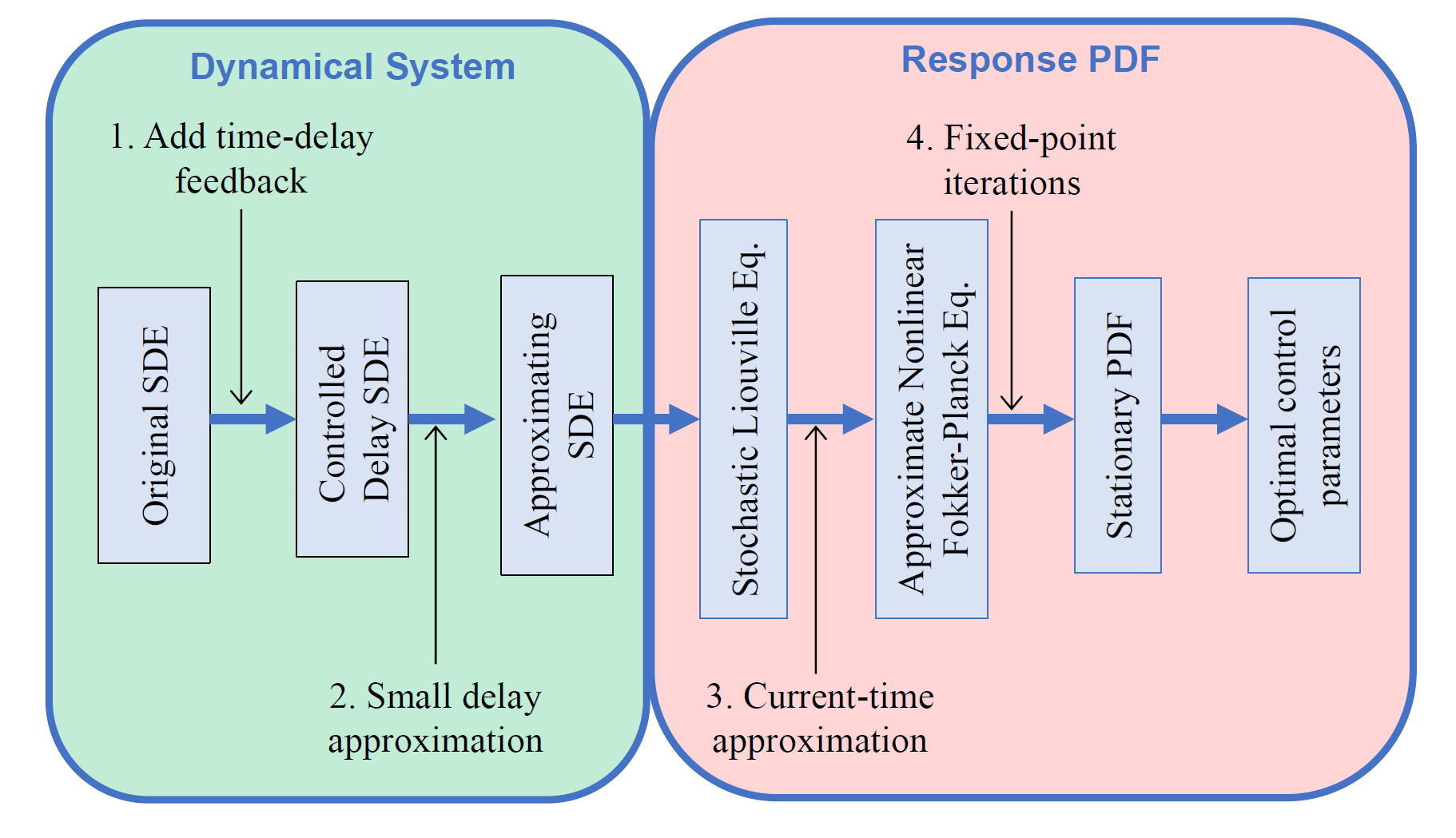}
	\caption{Summary of the program followed in this paper.}
	\label{fig:summary}
\end{figure*}

Figure~\ref{fig:summary} summarizes the program of this paper which we now briefly outline.
\begin{enumerate}
	\item Time-delay feedback control: The uncontrolled system is described by a stochastic differential equation (SDE) driven by colored noise. To mitigate the rare transitions, we add a time-delay feedback control to the SDE. As a result, 
	the controlled system is a stochastic \emph{delay} differential equation (SDDE). 
	\item Small delay approximation: Assuming that the control delay is relatively small, we 
	use the Taylor expansion of the control term to approximate this SDDE with an appropriate SDE. 
	\item Current-time approximation and nonlinear Fokker--Planck:
	The approximating SDE is still driven by colored noise and therefore its PDF evolution is described 
	by a non-closed, stochastic Liouville equation (see e.g. \cite[Sec. III.D]{Hanggi1995}). Our second approximation replaces the stochastic Liouville equation with a
	nonlinear Fokker--Planck equation in closed form.
	\item Fixed-point iterations for the stationary PDF: The stationary solution of the nonlinear Fokker--Planck equation is known for one-dimensional SDEs. However, this is an implicit solution as it depends on the response moment
	of the system, which itself depends on the stationary PDF. We devise a rapidly-converging iterative algorithm
	in order to estimate the stationary PDF and the response moment simultaneously.
\end{enumerate}

The above series of approximations allow for estimating the stationary response PDF in an inexpensive manner
which in turn enables us to determine the optimal control parameters parsimoniously. As we show with detailed numerical 
simulations, in spite of several approximations, the estimated PDFs agree remarkably well with the true PDFs 
obtained from Monte Carlo simulations.

One of the main contributions of the present work is the development of a new one-dimensional nonlinear Fokker--Planck equation whose stationary solution can be easily determined, and constitutes a fairly accurate approximation of the stationary response PDF of a scalar SDE under colored noise excitation. For additive noise excitation, a similar nonlinear Fokker--Planck equation was derived by Mamis et al.~\cite{Mamis2019d}; here we generalize the equation to the case of multiplicative 
stochastic excitations (Section~\ref{sec:Fokker--Planck}).

The remainder of this paper is organized as follows.
In Section~\ref{sec:prelim}, we describe the set-up of the problem and review some mathematical preliminaries.
In Section~\ref{sec:controlled_dynamical_systems}, we introduce the time-delay feedback control and discuss its 
effect on the system.
In Section~\ref{sec:control_param}, we derive the nonlinear Fokker--Planck equation, devise an iterative algorithm
to approximate its stationary solution, and discuss how it enables us to determine the optimal control parameters.
Section~\ref{sec:numerics} contains our numerical results. We discuss two examples: a stochastic system driven
by additive noise (Section~\ref{sec:additive}) and an SDE arising in optical lasers which is driven by multiplicative noise (Section~\ref{sec:multiplicative}).
Finally, we present our concluding remarks in Section~\ref{sec:conc}.

\section{Preliminaries and set-up}\label{sec:prelim}
\subsection{Uncontrolled stochastic dynamical systems}
Dynamical systems, driven by a potential and under random noise excitations, can be described by multidimensional stochastic differential equations (SDEs) of the form 
\begin{equation} \label{eq:Ndim_SDE}
\frac{\id\bm{X}(t)}{\id t}=-\nabla V(\bm{X}(t))+\bm{\s}(\bm{X}(t))\bm{\xi}(t), \ \ \bm{X}(t_0)=\bm{x}^0,
\end{equation}
where $\bm{X}(t)\in\mathbb R^n$ is the stochastic process of the system's state at time $t$, $V:\R^n\rightarrow\R$ is the potential function, $\bm{\xi}(t)\in\mathbb R^m$ is the noise excitation and $\bm{\s}(\bm{x})\in\mathbb R^{n\times m}$ is the noise intensity. If the matrix $\bm{\s}$ is a constant, independent of $\vc x$, the excitation is called additive; whereas in the general case of state-dependent function $\bm{\s}(\bm{x})$, the excitation is called multiplicative. 

It is often assumed that the noise $\bm{\xi}(t)$ is the standard multidimensional Gaussian white noise $\bm{\xi}^{\text{WN}}(t)$, with independent components, zero mean value and two-time autocorrelation  matrix
\begin{equation} \label{eq:cor_Ndim_WN}
\bm{C}_{\bm{\xi}}^{\text{WN}}(t_1,t_2)=\mathbb{E}\left[\bm{\xi}^{\text{WN}}(t_1)\left(\bm{\xi}^{\text{WN}}(t_2)\right)^T\right]=\bm{I}\,\delta(t_1-t_2),
\end{equation}
where $\mathbb{E}[\cdot]$ is the expected value, $T$ superscript denotes the matrix transpose, $\bm{I}$ is the $m\times m$ identity matrix and $\delta(t_1-t_2)$ is Dirac's delta function. White noise $\bm{\xi}^{\text{WN}}(t)$ is the formal time derivative of the standard Wiener process \cite{Pugachev2001, Horsthemke2006}.

Stochastic systems driven by white noise excitation have been studied extensively, with the development of both It{\=o} calculus \cite{Oksendal2003} for solving SDEs, and the formulation of corresponding Fokker--Planck equation \cite{Risken1996}. Fokker--Planck equation governs the system's response probability density function $p(\bm{x},t)$, defined so that probability $\mathbb{P}(\bm{X}(t)\in\mathcal{S})=\int_{{\mathcal{S}}}p(\bm{x},t)\id \bm{x}$ for any Lebesgue-measurable set $\mathcal{S}\subset\R^N$.
\subsection{Shaping filters}\label{sec:filters}
The description of environmental noises as white is not realistic. This can be easily seen by calculating the Fourier transform of its autocorrelation, which results in a constant, diagonal power spectrum matrix with infinite bandwidth (see Refs.~\cite[Sec. 3.2]{Horsthemke2006} and \cite[Sec. 1.4.2]{Gardiner2004}). A more realistic noise can be obtained by using shaping filters~\cite{Pugachev2001,Roberts2003}. Shaping filters are SDEs with white noise input, whose response is a colored noise, i.e. a smoothly-correlated stochastic process with a prescribed spectrum, or equivalently, a prescribed autocorrelation function. Following \cite[Sec. 5.4.2]{Pugachev2001}, we express a shaping filter as the general nonlinear SDE,
\begin{equation}\label{eq:filter}
\frac{\id\bm{\xi}(t)}{\id t}=\bm{\alpha}(\bm{\xi}(t))+\bm{\beta}(\bm{\xi}(t))\bm{\xi}^{\text{WN}}(t),
\end{equation}
where $\bm \alpha:\mathbb R^m\to \mathbb R^m$ and $\bm \beta:\R^m\to\R^{m\times m}$ are prescribed drift and diffusion functions, respectively.
The filters commonly used in applications are linear in the sense that the drift and diffusion are linear functions~\cite{Spanos1986a, Francescutto2004, Scruggs2013, Chai2015}. 

For instance, the \textit{Ornstein-Uhlenbeck} (OU) process is obtained when $\bm{\alpha}(\bm{\xi}(t))=-\bm{A}\bm{\xi}(t)$ and $\bm{\beta}(\bm{\xi}(t))=\bm{B}$, with $\bm{A}$, $\bm{B}$ being constant matrices. If all eigenvalues of drift matrix $\bm{A}$ have positive real parts, the stationary two-time autocorrelation function of OU process is the $m\times m$ matrix~\cite[ Sec. 4.4.6]{Gardiner2004}
\begin{align}\label{eq:OU_multi_cor}
\bm{C}_{\bm{\xi}}^{\text{OU}}(t_1,t_2)=\left\{
\begin{array}{ll}
      e^{-\bm{A}(t_1-t_2)}\bm{\Sigma}, & t_1>t_2, \\
      \bm{\Sigma}e^{-\bm{A}^T(t_2-t_1)} & t_1<t_2,\\
\end{array} 
\right. 
\end{align}
where $\bm{\Sigma}$ is the stationary OU covariance matrix, given by the algebraic Lyapunov equation
\begin{equation}\label{eq:OU_cov}
\bm{A\Sigma}^T+\bm{\Sigma A}=\bm{BB}^T.
\end{equation}
An important special case is the multidimensional OU noise with independent components, for which both drift and difussion matrices $\bm{A}$, $\bm{B}$ are equal to the diagonal matrix
\begin{equation}\label{eq:A_B}
\bm{A}=\bm{B}=\text{diag}\left[1/s_{cor}^{(1)},\ldots,1/s_{cor}^{(m)}\right],
\end{equation}
where $s_{cor}^{(\ell)}$ is the \textit{correlation time} of the $\ell$-th OU component. It is easy to see that, under Eq.~\eqref{eq:A_B}, OU autocorrelation~\eqref{eq:OU_multi_cor} simplifies to
\begin{align}
&\bm{C}_{\bm{\xi}}^{\text{OU}}(t_1,t_2)=\nonumber\\ &\text{diag}\left[\frac{1}{2s_{cor}^{(1)}}\exp\left(-\frac{|t_1-t_2|}{s_{cor}^{(1)}}\right),\ldots,\frac{1}{2s_{cor}^{(m)}}\exp\left(-\frac{|t_1-t_2|}{s_{cor}^{(m)}}\right)\right]. \label{eq:OU_multi_cor2}
\end{align} 
Henceforth, OU noise with autocorrelation~\eqref{eq:OU_multi_cor2} will be called the \textit{standard multidimensional OU noise}. Since $\left(1/\varepsilon\right)\exp\left(-|t_1-t_2|/\varepsilon\right)$ is a nascent delta function \cite[Sec. 6.6]{Toral2014}, standard OU autocorrelation~\eqref{eq:OU_multi_cor2} tends to white noise autocorrelation~\eqref{eq:cor_Ndim_WN} for $\max\left\{s_{cor}^{(1)},\ldots,s_{cor}^{(m)}\right\}\rightarrow 0$. Thus, standard OU colored noise tends to white noise when all its correlation times go to zero.

Note that, by considering the system consisting of the original SDE~\eqref{eq:Ndim_SDE} and shaping filter SDE~\eqref{eq:filter}, we obtain a multidimensional augmented SDE. 
Therefore, the original SDE excited by colored noise is replaced by a higher dimensional SDE excited by white noise. However, unlike the original SDE~\eqref{eq:Ndim_SDE}, the drift term
$[-\nabla V(\bm{X}), \pmb\alpha(\pmb\xi)]^T$
in the augmented SDE cannot always be expressed as the gradient of a potential function. 

\section{Controlled Stochastic Dynamical Systems}\label{sec:controlled_dynamical_systems}
In the absence of noise, the minima of the potential $V$ are stable equilibria of system~\eqref{eq:Ndim_SDE}. 
We assume that one of these equilibria is desirable and denote its position with $\vc x_a$.
All other equilibria are undesirable or `bad' and are denoted by $\vc x_b$. Ideally, we would like the system to 
evolve in the vicinity of the desirable equilibrium $\vc x_a$.
However, the presence of noise $\bxi(t)$ enables rare transitions away from the desirable equilibrium towards
an undesirable equilibrium $\vc x_b$. Since these transitions can have catastrophic consequences,
we would like to design a simple control strategy that mitigates transitions away from the desirable equilibrium $\vc x_a$.

In order to be able to suppress rare transitions away from $x_a$, we add the \textit{linear, time-delay feedback control term} $-a\left(\bm{X}(t-\tau)-\hat{\bm{x}}\right)$ to the original SDE~\eqref{eq:Ndim_SDE}. This results in the controlled stochastic delayed differential equation (SDDE) \cite{Farazmand2020}
\begin{equation} \label{eq:Ndim_SDDE}
\frac{\id\bm{X}(t)}{\id t}=-\nabla V(\bm{X}(t))-a(\bm{X}(t-\tau)-\hat{\bm{x}})+\bm{\s}(\bm{X}(t))\bm{\xi}(t),
\end{equation}
supplemented by the appropriate initial condition; $\bm{X}(t)=\bm{x}^0(t)$, for $t\in[-\tau,0]$. In the linear control term, $a$ is the control gain, $\tau>0$ is the time delay, and $\hat{\bm{x}}$ is the shift chosen in order to suppress transitions away from the desirable equilibrium $\bm{x}_a$ in the controlled system response. 

We chose this linear delayed feedback as control term because it is easy to implement, and does not give rise to additional nonlinearities in the dynamical system. This control strategy has been widely employed in the stabilization of deterministic systems (see e.g.~\cite{Pyragas1995, Suresh2018}). Also, there is a practical reason why we use the delayed response $\bm{X}(t-\tau)$ in the control term, instead of $\bm{X}(t)$; the delay $\tau$ models the inevitable lag between state measurements and control actuation.

In Ref.~\cite{Farazmand2020}, where the control of SDEs under white noise was studied, $\hat{\bm{x}}$ was chosen equal to the desirable equilibrium $\bm{x}_a$. However, in the colored excitation case, the choice of $\hat{\bm{x}}$ is not so straightforward. This is due to the appearnce of peak drift phenomenon, which will be discussed at length in sections \ref{sec:additive} and \ref{sec:multiplicative}. This means that, in order to be able to suppress the peak drift, we choose $\hat{\bm{x}}$ to be in the vicinity of $\bm{x}_a$, but not exactly equal to $\bm{x}_a$. For this, in our analysis, we do not assume that $\hat{\bm{x}}$ and $\bm{x}_a$ coincide. 

By confining ourselves to the small time delay regime,  $0<\tau\ll 1$, the delayed term $\bm{X}(t-\tau)$ in SDDE~\eqref{eq:Ndim_SDDE} can be approximated by a linear Taylor expansion with respect to $\tau$:
\begin{equation}
\bm{X}(t-\tau)= \bm{X}(t)-\tau\dot{\bm{X}}(t)+\mathcal{O}(\tau^2). \label{eq:tau_taylor}
\end{equation}
By neglecting $\mathcal{O}(\tau^2)$ terms, this gives rise to the approximating SDE
\begin{equation} \label{eq:Ndim_approx_SDE}
(1-a\tau)\frac{\id\bm{X}(t)}{\id t}=-\nabla\tilde{V}(\bm{X}(t))+\bm{\s}(\bm{X}(t))\bm{\xi}(t),
\end{equation}
where $\tilde{V}(\bm{x})$ is the \textit{effective potential}, defined as
\begin{equation}\label{eq:eff_potential}
\tilde{V}(\bm{x})=V(\bm{x})+\frac{a}{2}(\bm{x}-\hat{\bm{x}})^2.
\end{equation}
Therefore, the control modifies the effective potential of the system. The following theorem shows that 
the control also effectively modifies the filter producing the colored noise excitation. 
\begin{thm}[The rescaled approximating SDE under general colored noise] \label{th_rescaling}
Consider the system consisting of the approximating SDE~\eqref{eq:Ndim_approx_SDE} and noise filter SDE~\eqref{eq:filter}. Defining the rescaled time 
$s=t/(1-a\tau)$, the SDE system can be expressed equivalently as
\begin{align} 
\frac{\id \bm{X}(s)}{\id s}&=-\nabla\tilde{V}(\bm{X}(s))+\bm{\s}(\bm{X}(s))\bm{\xi}(s) \label{eq:Ndim_resc_SDE}\\
\frac{\id \bm{\xi}(s)}{\id s}&=(1-a\tau)\bm{\alpha}(\bm{\xi}(s))+\sqrt{1-a\tau}\bm{\beta}(\bm{\xi}(s))\bm{\xi}^{\text{WN}}(s), \label{eq:Ndim_resc_filter}
\end{align}
where $\bm{\xi}^{\text{WN}}(s)$ is a standard white noise.
\end{thm}
\begin{proof}
See Appendix \ref{A_rescaling}.
\end{proof}

The implications of this theorem become clear if we assume $\bxi(t)$ is the standard OU noise excitation as discussed in Section~\ref{sec:filters}.
Therefore, in the following corollary, we use the above general result in order to obtain the rescaled approximating SDE for the case of standard OU noise excitation.  
\begin{cor}[The rescaled approximating SDE under standard OU noise] \label{cor:rescaled_OU} 
For the case of standard OU noise excitation with autocorrelation~\eqref{eq:OU_multi_cor2}, approximating SDE~\eqref{eq:Ndim_approx_SDE} is equivalent to the rescaled SDE:
\begin{equation}\label{eq:Ndim_SDE_OU_aug}
\frac{\id \bm{X}(s)}{\id s}=-\nabla\tilde{V}(\bm{X}(s))+\tilde{\bm{\s}}(\bm{X}(s))\tilde{\bm{\xi}}^{\text{OU}}(s),
\end{equation}
where \textit{effective noise intensity} is given by
\begin{equation}\label{eq:tilde_sigma}
\tilde{\bm{\s}}(\bm{x})=\frac{\bm{\s}(\bm{x})}{\sqrt{1-a\tau}},
\end{equation} 
and $\tilde{\bm{\xi}}^{\text{OU}}(s)$ is the rescaled standard OU noise with correlation times
\begin{equation}\label{eq:tilde_scor}
\tilde{s}_{cor}^{(\ell)}=\frac{s_{cor}^{(\ell)}}{1-a\tau}, \ \ \ \ell=1,\ldots,m.
\end{equation}
 \end{cor}
 \begin{proof}
This follows from a direct application of Theorem \ref{th_rescaling}; see Appendix \ref{A:cor} for details.
\end{proof}
\begin{rem}
Thus, the effects of control on an SDE excited by standard OU noise are the following:
\begin{enumerate}
\item The control term modifies the potential $V(\bm{x})$ of the original SDE to effective potential $\tilde{V}(\bm{x})$, given by Eq.~\eqref{eq:eff_potential}. Since value $\hat{\bm{x}}$ is equal or in the vicinity of desirable equilibrium $\bm{x}_a$, the term $(a/2)(\bm{x}-\hat{\bm{x}})^2$ deepens the potential well around $\bm{x}_a$, and therefore hinders escapes away from it.
\item The control also intensifies noise, since the noise intensity is multiplied by the factor $1/\sqrt{1-a\tau}$ (see Eq.~\eqref{eq:tilde_sigma}). This effect is antagonistic to the stabilization around $\bm{x}_a$, since higher noise levels can result in transitions away from the desirable equilibrium.
This increase in effective noise intensity was first reported by Guillouzic et al.~\cite{Longtin1999} and later redicsovered by Farazmand~\cite{Farazmand2020} using a different approach.
\item Finally, the control increases correlation times of standard OU noise by a factor of $1/(1-a\tau)$, see Eq.~\eqref{eq:tilde_scor}. This increase in the noise correlation time renders the white noise approximation even more inapplicable to the controlled SDE case.
\end{enumerate}
\end{rem}
\begin{rem}[The limiting case of white noise excitation]\label{rem:white_noise}
As we have mentioned in Sec. \ref{sec:filters}, standard OU noise tends to white noise when its correlation times go to zero. Thus, for $\max\left\{s_{cor}^{(1)},\ldots,s_{cor}^{(m)}\right\}\rightarrow 0$, rescaled approximating SDE~\eqref{eq:Ndim_SDE_OU_aug} results in
\begin{equation}\label{eq:Ndim_SDE_WN_aug}
\frac{\id \bm{X}(s)}{\id s}=-\nabla\tilde{V}(\bm{X}(s))+\tilde{\bm{\s}}(\bm{X}(s))\tilde{\bm{\xi}}^{\text{WN}}(s),
\end{equation}
which is the rescaled approximating SDE for the white noise excitation case obtaned in Ref. \cite{Farazmand2020}.
\end{rem}

While the above results describe the effects of feedback delay control on multidimensional stochastic dynamical systems, we shall focus, for the rest of the present work, on the scalar case. The list of equations and parameters for the control of the scalar SDE is summarized as
\begin{itemize}
\item Uncontrolled SDE: \begin{equation} 
\frac{\id X(t)}{\id t}=-V'(X(t))+\s(X(t))\xi(t). \label{eq:SDE}
\end{equation}
\item Scalar OU noise excitation: \begin{equation}
C^{\text{OU}}_{\xi}(t_1,t_2)=\frac{1}{2s_{cor}}\exp\left(-\frac{|t_1-t_2|}{s_{cor}}\right).
\label{eq:corr_OU}
\end{equation}
\item Controlled SDDE: \begin{equation}
\frac{\id X(t)}{\id t}=-V'(X(t))-a\left(X(t-\tau)-\hat{x}\right)+\s(X(t))\xi^{\text{OU}}(t).
\label{eq:SDDE}
\end{equation}
\item Approximating SDE: \begin{equation}
(1-a\tau)\frac{\id X(t)}{\id t}=-\tilde{V}'(X(t))+\s(X(t))\xi^{\text{OU}}(t),
\label{eq:approx_SDE}
\end{equation}
with $\tilde{V}(x)=V(x)+(a/2)(x-\hat{x})^2$.
\item Rescaled approximating SDE: \begin{equation}
\frac{\id X(s)}{\id s}=-\tilde{V}'(X(s))+\tilde{\s}(X(s))\tilde{\xi}^{\text{OU}}(s), \label{eq:rescaled_SDE}
\end{equation}
with $s=t/(1-a\tau)$, $\tilde{\s}(x)=\s(x)/\sqrt{1-a\tau}$, and $\tilde{s}_{cor}=s_{cor}/(1-a\tau)$.
\end{itemize}
The reason for choosing a scalar SDE is that, in this case, the stationary PDF is readily available from the solution of either the classical or the nonlinear Fokker--Planck equations, as we will see in the Section~\ref{sec:control_param}. Thus, for scalar SDE~\eqref{eq:SDE}, we can present as well as show the validity of our methodology for mitigating rare events, without the need for a numerical solver of the nonlinear Fokker--Planck equation.

\section{Determining the control parameters}\label{sec:control_param}
	The optimal values of the control parameters $(a,\tau,\hat x)$, can be determined 
	by direct Monte Carlo simulations. To avoid such computationally expensive simulations, we determine the
	optimal control parameters by studying the stationary response PDF of the system. In Section~\ref{sec:Fokker--Planck}, we 
	first discuss the nonlinear Fokker--Planck equations which approximate the PDF evolution
	of an SDE excited by colored noise. Subsequently, in Section~\ref{sec:stationary_pdf}, we devise an iterative algorithm for
	obtaining the stationary PDF of the nonlinear Fokker--Planck equation.
\subsection{Nonlinear Fokker--Planck equation}\label{sec:Fokker--Planck}
It is well-established, see e.g. \cite[Chapter 5]{Gardiner2004}, that the evolution of response PDF $p(x,t)$ of SDE~\eqref{eq:SDE} under white noise excitation, $\xi(t)=\xi^{\text{WN}}(t)$, is governed by the classical one-dimensional Fokker--Planck equation,
\begin{align}
\frac{\d p(x,t)}{\d t}=&\frac{\d}{\d x}\left\{\left[V'(x)-\frac{\varpi}{2}\s'(x)\s(x)\right]p(x,t)\right\}+\nonumber\\ &+\frac{1}{2}\frac{\d^2}{\d x^2}\left[\s^2(x)p(x,t)\right]. \label{eq:Fokker--Planck}
\end{align}
\begin{rem}[The Wong--Zakai correction] \label{rem:WZ}
In Fokker--Planck Eq.~\eqref{eq:Fokker--Planck}, the drift coefficient is augmented with the term $(1/2)\s'(x)\s(x)$, which is the Wong--Zakai correction, modeling the difference between the It{\=o} ($\varpi=0$) and Stratonovich ($\varpi=1$) interpretations of SDEs under white noise \cite{Sun2006, Mamis2016e}. Since Wong--Zakai correction depends on $\s'(x)$, It{\=o} and Stratonovich Fokker--Planck equations coincide in the case of additive white noise excitation, where $\s$ is constant.
\end{rem}

This convenient description, via a single one-dimensional partial differential equation of drift-diffusion type, is not readily applicable to SDEs under colored noise excitation. As discussed in \ref{sec:filters}, using shaping filters, one can still express the SDE driven by colored noise as an equivalent multidimensional augmented SDE driven by white noise.
However, the drift term in the augmented SDE is no longer the gradient of a potential, and consequently, the closed-form solution of the resulting Fokker--Planck equation is generally unknown~\cite{Masud2005, Chen2017, Xu2020}. 

In the present work, we take an alternative path; we use an approximate Fokker--Planck-like equation corresponding to SDE~\eqref{eq:SDE} driven by colored noise, without resorting to shaping filters. This approximate Fokker--Planck equation is easily solvable, and allows us to determine the appropriate 
control parameters without resorting to computationally expensive Monte Carlo simulations.

There exists an extensive body of work devoted to deriving approximate Fokker--Planck equations for SDEs excited by colored noise
(see, e.g.,~\cite{Fox1987, Faetti1988, Peacock-Lopez1988, Hanggi1989, Hanggi1995, Ridolfi2011, Bianucci2020a}, to mention a few main references). 
An approximate Fokker--Planck equation that is still used in applications 
is the one proposed by H{\"a}nggi and his team~\cite{Hanggi1985,Hanggi1995}. 
Recently, Mamis et al.~\cite{Mamis2019d} generalized H{\"a}nggi's equation
by incorporating higher-order corrections and therefore obtaining more accurate response probability densities.
Here, we generalize the derivation of Ref.~\cite{Mamis2019d}  to include the more general case of multiplicative excitation.

The approximate, Fokker--Planck-like, response PDF evolution equations that correspond to SDE~\eqref{eq:SDE}, with differentiable and non-vanishing ($\s(x)\neq 0$) noise intensity and twice continuously differentiable potential $V(x)$, and under Gaussian colored noise excitation $\xi(t)$ with a general, non-singular autocorrelation function $C_{\xi}(t,s)$, read
\begin{align}
\frac{\d p(x,t)}{\d t}=&\frac{\d}{\d x}\left\{\left[V'(x)-\s'(x)\s(x)A_M(x,t;p)\right]p(x,t)\right\}+\nonumber \\&+\frac{\d^2}{\d x^2}\left[\s^2(x)A_M(x,t;p)p(x,t)\right]. \label{eq:nFokker--Planck}
\end{align} 
The coefficient $A_M$ is defined, for $M=0\text{ or }2$, as
\begin{equation}
A_M(x,t;p)=\sum_{m=0}^M\frac{D_m(t;p)}{m!}\big\{\zeta(x)-\mathbb{E}\left[\zeta(X(t))\right]\big\}^m, \label{eq:diff}
\end{equation}
where
\begin{equation}
 \zeta(x)=-\s(x)\left(\frac{V'(x)}{\s(x)}\right)'\label{eq:zeta1}
 \end{equation}
 and 
\begin{equation}
D_m(t;p)=\int_{t_0}^tC_{\xi}(t,t_1)\exp\int_{t_1}^t\mathbb{E}\left[\zeta(X(u))\right]\id u(t-t_1)^m\id t_1. \label{eq:D}
\end{equation}
The derivation of Eq.~\eqref{eq:nFokker--Planck} from the stochastic Liouville equation~\cite[Sec. III.D]{Hanggi1995} is performed in Appendix \ref{A:1}. For $M=0$, Eq.~\eqref{eq:nFokker--Planck} results in the usual H{\"a}nggi equation, while, for $M=2$, a new evolution equation is obtained. As shown in Section~\ref{sec:numerics}, the case with $M=2$ consistently renders more accurate results than H{\"a}nggi's $M=0$ approximation. 

The main difference between PDF evolution equations given by~\eqref{eq:nFokker--Planck} and the classical Fokker--Planck Eq.~\eqref{eq:Fokker--Planck} is that, in Eq.~\eqref{eq:nFokker--Planck} , coefficient $A_M$ depends not only on state variable $x$ and time $t$, but also on the unknown response PDF $p$. More specifically, $A_M$ depends on the time history of the response moment
\begin{equation}
\mathbb{E}\left[\zeta(X(t))\right]=\int_{\R}\zeta(x)p(x,t)\id x, \label{eq:zeta_moment}
\end{equation}
through Eqs.~\eqref{eq:diff} and~\eqref{eq:D}. Such equations are commonly referred to as \textit{nonlinear Fokker--Planck equations} \cite{Frank2005}. 

Considering the Gaussian colored excitation $\xi(t)$ of the standard OU noise with autocorrelation~\eqref{eq:corr_OU},
the coefficients $D_m$ read
\begin{align}
D_m^{\text{OU}}&(t;p)=\frac{1}{2s_{cor}}\int_{t_0}^t\exp\left(-\frac{t-t_1}{s_{cor}}\right)\times \nonumber\\ &\times\exp\left(\int_{t_1}^t\mathbb{E}\left[\zeta(X(u))\right]\id u\right)(t-t_1)^m\id t_1. \label{eq:D_OU}
\end{align}
From this point forward, $\xi(t)$ will denote OU noise, unless it is explicitly stated otherwise. 

\begin{rem}[Compatibility with time rescaling of corollary \ref{cor:rescaled_OU}]
	Straightforward algebraic manipulations show that, by rescaling time $s=t/(1-a\tau)$ in the nonlinear Fokker--Planck Eq.~\eqref{eq:nFokker--Planck} corresponding to approximating controlled SDE~\eqref{eq:approx_SDE} under OU noise, the nonlinear Fokker--Planck equation corresponding to rescaled SDE~\eqref{eq:rescaled_SDE} is obtained. 
\end{rem}

\begin{cor}[Compatibility with classical Fokker--Planck equation]\label{rem:compatibility}
In the limiting white noise case $s_{cor}\rightarrow0$, the nonlinear Fokker--Planck Eq.~\eqref{eq:nFokker--Planck} coincides with Stratonovich's Fokker--Planck Eq.~\eqref{eq:Fokker--Planck} with $\varpi=1$.
\end{cor}
\begin{proof}
Applying a change of variable, Eq.~\eqref{eq:D_OU} is expressed equivalently as
\begin{align}
D_m^{\text{OU}}(t;p)&=\frac{1}{2}\int_0^{\frac{t-t_0}{s_{cor}}}\exp\left(-v\right)q_m(vs_{cor})\id v\nonumber\\&\equiv\frac{1}{2}\int_0^{+\infty}\exp\left(-v\right)q_m(vs_{cor})I_{\left(0,\frac{t-t_0}{s_{cor}}\right)}(v)\id v, \label{eq:D_OU3}
\end{align}
with $q_m(s):=\exp\left(\int_0^s\mathbb{E}\left[\zeta(X(u))\right]\id u\right)s^m$, and $I_{\left(0,\frac{t-t_0}{s_{cor}}\right)}(v)$ being the indicator function that takes the value 1 for $v\in\left(0,\frac{t-t_0}{s_{cor}}\right)$, and 0 otherwise. Applying the dominated convergence theorem to~\eqref{eq:D_OU3}, we obtain
\begin{equation}
\lim_{s_{cor}\rightarrow0}D_m^{\text{OU}}(t;p)=\frac{q_m(0)}{2}\int_0^{+\infty}\exp\left(-v\right)\id v=\frac{q_m(0)}{2}. \label{eq:D_OU4}
\end{equation}
Since $q_0=1$ and $q_1=q_2=0$, we obtain $D_0=1/2$, $D_1(t)=D_2(t)=0$ after taking the limit $s_{cor}\rightarrow 0$, resulting in $A_M=1/2$. 
\end{proof}

Corollary~\ref{rem:compatibility} is compatible with previous results; that is, if the white noise excitation is approximated by a colored noise with a very short but nonzero correlation time, Stratonovich’s interpretation of the SDE should be applied (see, e.g.,~\cite[page 128]{Ottinger1996} and \cite[page 216]{Horsthemke2006}). 

\subsection{Stationary distribution}\label{sec:stationary_pdf}
In this section, we discuss on the stationary response PDF $p_0(x)$, which is invariant in time 
and describes the long-term evolution of the SDE response process.
This stationary PDF is immensely helpful for determining the optimal control parameters, 
without resorting to computationally expensive Monte Carlo simulations.

The following lemma determines the stationary form of the nonlinear Fokker--Planck Eq.~\eqref{eq:nFokker--Planck} for OU excitation.
\begin{lem}\label{lem:stationary_nFokker--Planck}
Assume that response moment $R:=\mathbb{E}\left[\zeta(X(t))\right]$ attains a finite time-independent value, satisfying the condition
\begin{equation}
R<\frac{1}{s_{cor}}. \label{eq:stationary_condition}
\end{equation}
Then, the stationary nonlinear Fokker--Planck Eq.~\eqref{eq:nFokker--Planck} for SDE~\eqref{eq:SDE} under OU excitation reads
\begin{align}
\frac{\id}{\id x}&\left\{\left[V'(x)-\s'(x)\s(x)A_M(x,R)\right]p_0(x)\right\}\nonumber\\ &+\frac{\id^2}{\id x^2}\left[\s^2(x)A_M(x,R)p_0(x)\right]=0, \label{eq:stationary_nFokker--Planck}
\end{align}
where the stationary coefficient $A_M(x,R):=\lim_{t\rightarrow\infty}A_M(x,t;p)$ is given by
\begin{equation}
A_M(x,R)=\frac{1}{2}\sum_{m=0}^M\frac{\left[s_{cor}\left(\zeta(x)-R\right)\right]^m}{\left(1-s_{cor}R\right)^{m+1}}. \label{eq:stationary_A_M}
\end{equation}  
\end{lem}
\begin{proof}
It is easy to see that, under condition~\eqref{eq:stationary_condition}, $D_m^{\text{OU}}$ stationary values are finite and can be calculated from their definition~\eqref{eq:D} as $t\rightarrow\infty$. By substituting this result into Eq.~\eqref{eq:diff}, relation~\eqref{eq:stationary_A_M} is obtained.
\end{proof}
\begin{rem}[Positivity of diffusion coefficient]\label{rem:positivity_diff}
For~\eqref{eq:stationary_nFokker--Planck} to be a valid stationary Fokker--Planck-like equation, its diffusion coefficient $\s^2(x)A_M(x,R)$ should be positive. For H{\"a}nggi's stationary equation, $M=0$, we have $A_0(R)=1/[2(1-s_{cor}R)]$ which is always positive under condition~\eqref{eq:stationary_condition}. For $M=2$, $A_2(x,R)$ can be written equivalently as
\begin{align}
A_2(x,R)=&\frac{1}{2(1-s_{cor}R)^3}\times\nonumber\\ \times[&3s_{cor}^2R^2-3s_{cor}(s_{cor}\zeta(x)+1)R+\nonumber\\ &+s^2_{cor}\zeta^2(x)+s_{cor}\zeta(x)+1]. \label{eq:A2}
\end{align}
On the right-hand side of Eq.~\eqref{eq:A2}, the fraction is always positive under condition~\eqref{eq:stationary_condition}. The other factor is identified as a quadratic polynomial with respect to $R$. By calculating its discriminant, $\Delta=-3s_{cor}^2(s_{cor}\zeta(x)-1)^2\leq 0$, we see that this polynomial always has the sign of coefficient $3s_{cor}^2$. Therefore, the diffusion coefficient is always positive. Note that positivity of the diffusion coefficient is not guaranteed in other approximate stationary PDF equations (see, for example, the Fokker--Planck-like equation derived under the small correlation time approximation in Refs.~\cite{Sancho1982} and~\cite[Sec. 8.6]{Horsthemke2006}.)
\end{rem}
\begin{cor}[Solution to stationary nonlinear Fokker--Planck equations]\label{cor:p_0}
The solution to the stationary Fokker--Planck-like equation~\eqref{eq:stationary_nFokker--Planck} is given by
\begin{align}
p_0(x,R)=& \frac{C(R)}{|\s(x)|A_M(x,R)}\times\nonumber\\
& \times\exp\left(-\int^x\frac{V'(y)}{\s^2(y)A_M(y,R)}\id y\right), \label{eq:p_0}
\end{align}
where $\int^x\id y$ denotes the antiderivative and $C(R)$ is the normalization factor, so that $\int_{\R}p_0(x,R)\id x=1$.
\end{cor}
\begin{proof}
See \cite[Sec. 5.2.2]{Gardiner2004} and \cite[Sec. 6.1]{Horsthemke2006}.
\end{proof}
We note that equation~\eqref{eq:p_0} is an implicit closed-form solution for the stationary nonlinear Fokker--Planck equation, since $p_0(x,R)$ depends on the response moment $R$ which remains to be determined. Determining $R$ in turn requires the knowledge of the stationary PDF $p_0$.
This is in contrast with the classical Fokker--Planck equation, where $A_M=1/2$ is independent of $R$, and therefore equation~\eqref{eq:p_0} constitutes its explicit closed-form solution. Thus, the dependence of response PDF $p_0$ on a response moment $R$ is a property arising from the colored excitation.

Nonetheless, Definition~\ref{defn:iter_scheme} below establishes an iterative scheme for calculating the response moment $R$ and therefore
the stationary PDF $P_0$.
\begin{defn}[The self-consistency equation for $R$]\label{defn:iter_scheme}
Using the definition of the response moment,
\begin{equation}
R=\int_{\R}\zeta(x)p_0(x,R)\id x, \label{eq:R_def}
\end{equation}
and the expression~\eqref{eq:p_0} for $p_0(x,R)$, we obtain the \textit{self-consistency equation} \cite{Frank2005},
\begin{equation}
R=\mathcal{I}(R), \label{eq:consistency}
\end{equation}
where
\begin{equation}\label{eq:I(R)}
\mathcal{I}(R)=\frac{\int_{\R}\frac{\zeta(x)}{|\s(x)|A_M(x,R)}\exp\left(-\int^x\frac{V'(y)}{\s^2(y)A_M(y,R)}\id y\right)\id x}{\int_{\R}\frac{1}{|\s(x)|A_M(x,R)}\exp\left(-\int^x\frac{V'(y)}{\s^2(y)A_M(y,R)}\id y\right)\id x}.
\end{equation}
\end{defn} 

Therefore, the correct value of the response moment $R$ is a fixed point of the map $\mathcal I:\R\to\R$.
If the map $\mathcal I$ is a contraction, the moment $R$ can be calculated through the iterative scheme, 
\begin{equation}
R_{n+1}=\mathcal{I}(R_n),\quad n=0,1,2,\cdots,
\end{equation}
so that $R=\lim_{n\to\infty}R_n$. We calculate an initial estimation $R_0$ from the explicit solution~\eqref{eq:p_0} for $A_M=1/2$ of the respective classical Fokker--Planck equation. Our iteration scheme is summarized in Algorithm \ref{algorithm}.
As seen in the numerical examples studied in Section~\ref{sec:numerics}, this iteration scheme is rapidly convergent. For instance, for an error tolerance $\varepsilon_{\text{tol}}=10^{-4}$, the scheme converges within 4 iterations on average.

\begin{algorithm}[H]
\caption{Iteration scheme for $R$} \label{algorithm}
\begin{enumerate}[leftmargin=*]
\item Choose a small tolerance $\varepsilon_{\text{tol}}$.
\item Calculate initial $R$ by substituting solution~\eqref{eq:p_0} with $A_M=\frac{1}{2}$ in response moment definition~\eqref{eq:R_def}
\item Determine $A_M(x)$ by substituting initial $R$ in Eq.~\eqref{eq:stationary_A_M} 
\item Calculate the update $R_{\text{upd}}=\mathcal I(R)$ using self-consistency Eq.~\eqref{eq:consistency}
\end{enumerate}
\textbf{while} $|R-R_{\text{upd}}|>\varepsilon_{\text{tol}}$ 
\begin{enumerate}[leftmargin=*, resume]
\item Set $R=R_{\text{upd}}$
\item Update $A_M(x)$ by substituting $R$ in Eq.~\eqref{eq:stationary_A_M} 
\item Calculate the next update $R_{\text{upd}}=\mathcal I(R)$ using self-consistency Eq.~\eqref{eq:consistency}
\end{enumerate}
\textbf{end while}
\begin{enumerate}[leftmargin=*, resume]
\item Determine stationary PDF $ p_0(x)$ using solution~\eqref{eq:p_0} and $R=R_{\text{upd}}$.
\end{enumerate}
\end{algorithm}

For a given set of control parameters $(a,\tau,\hat x)$, the stationary PDF of the approximating controlled SDE~\eqref{eq:rescaled_SDE}
can be readily computed through Algorithm \ref{algorithm}. Examining the response PDF $p_0$ determines whether the control has sufficiently
suppressed transitions away from the desirable equilibrium $x_a$. The optimal control parameters should result in a controlled response PDF that is
unimodal and concentrated around $x_a$. Since, as opposed to Monte Carlo simulations, 
computing the response PDF by Algorithm  \ref{algorithm} is computationally inexpensive, the control 
parameter space can be swept to determine the optimal control parameters.

\section{Numerical Results}\label{sec:numerics}
In this section, we present two numerical examples demonstrating the efficacy of the 
proposed time-delay feedback control. In Section~\ref{sec:additive}, we present the results for
a benchmark double-well potential driven by additive colored noise.
Section~\ref{sec:multiplicative} deals with a reduced-order model of an
optical laser which is driven by a multiplicative colored noise.

In this section, approximate stationary PDFs given by~\eqref{eq:p_0}, for $M=0$ (H{\"a}nggi's approximation) and for $M=2$ (our approximation), are compared to PDFs obtained from direct Monte Carlo (MC) simulations. To obtain the Monte Carlo results, the uncontrolled SDE and the corresponding controlled SDDE are augmented with the filter SDE for scalar OU noise, resulting in a two-dimensional SDE or SDDE driven by white noise (see the discussion on augmented systems in Section~\ref{sec:filters}). Trajectories of the resulting augmented systems are generated using the predictor-corrector scheme proposed by Cao et al.~\cite{Cao2015}. 
For the construction of each PDF from Monte Carlo simulations, $10^6$ realizations of the respective stochastic equation are used. 

\subsection{An additively excited bistable system} \label{sec:additive}
\begin{figure*}
		\subfigure[]{\label{fig:1a}\includegraphics[width=0.49\textwidth]{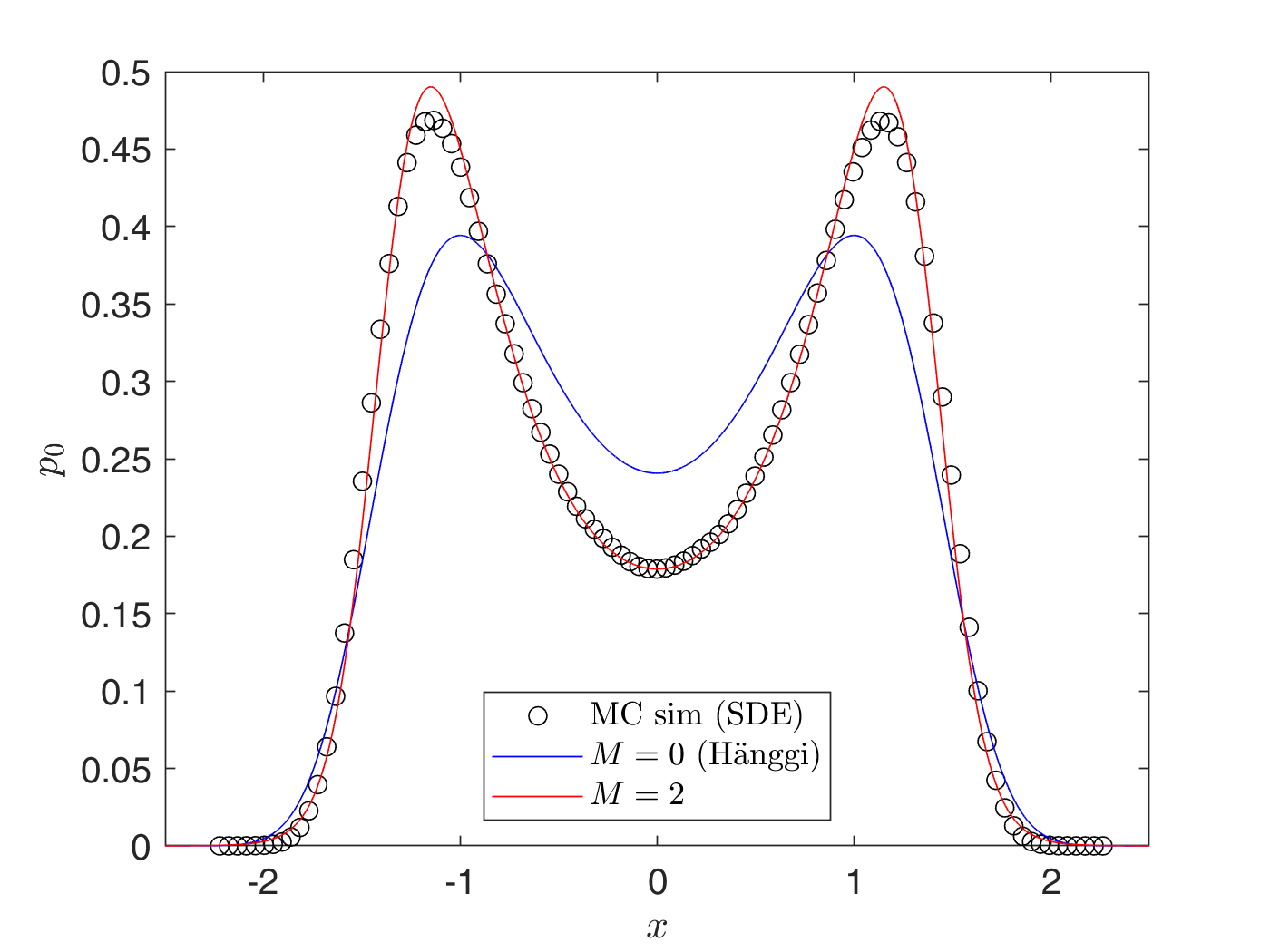}}
		\subfigure[]{\label{fig:1b}\includegraphics[width=0.49\textwidth]{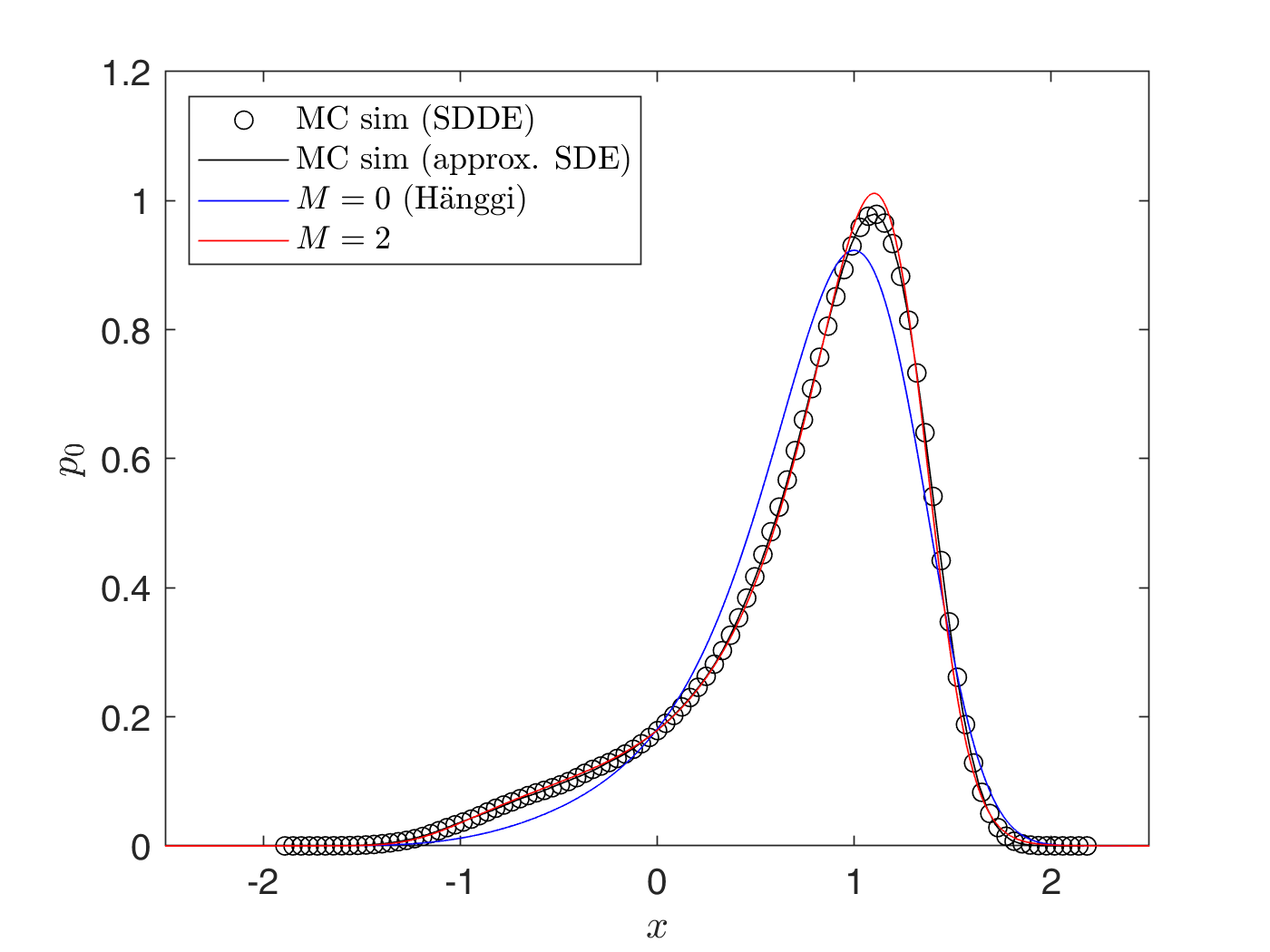}}
		\subfigure[]{\label{fig:1c}\includegraphics[width=0.49\textwidth]{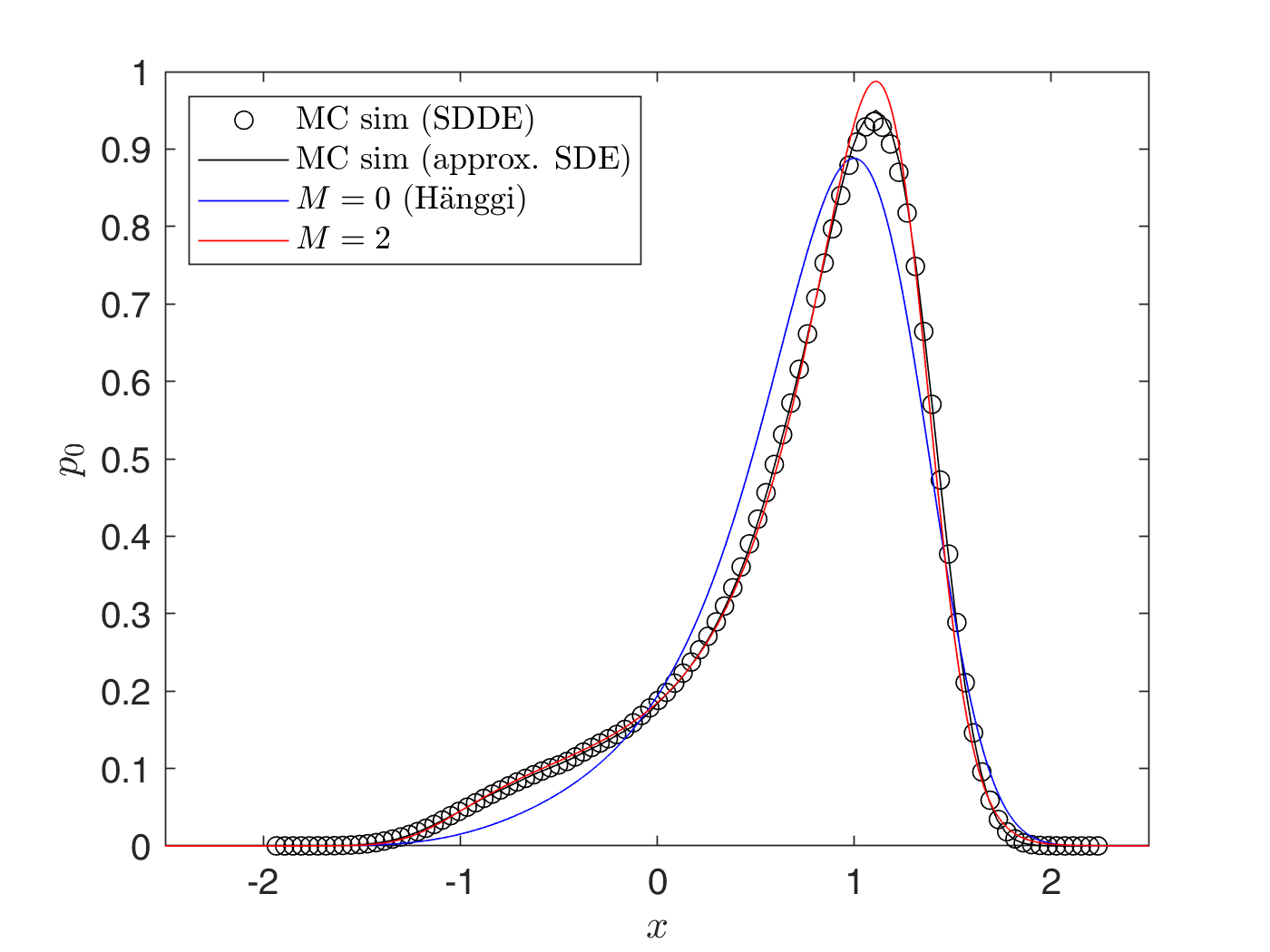}}
		\subfigure[]{\label{fig:1d}\includegraphics[width=0.49\textwidth]{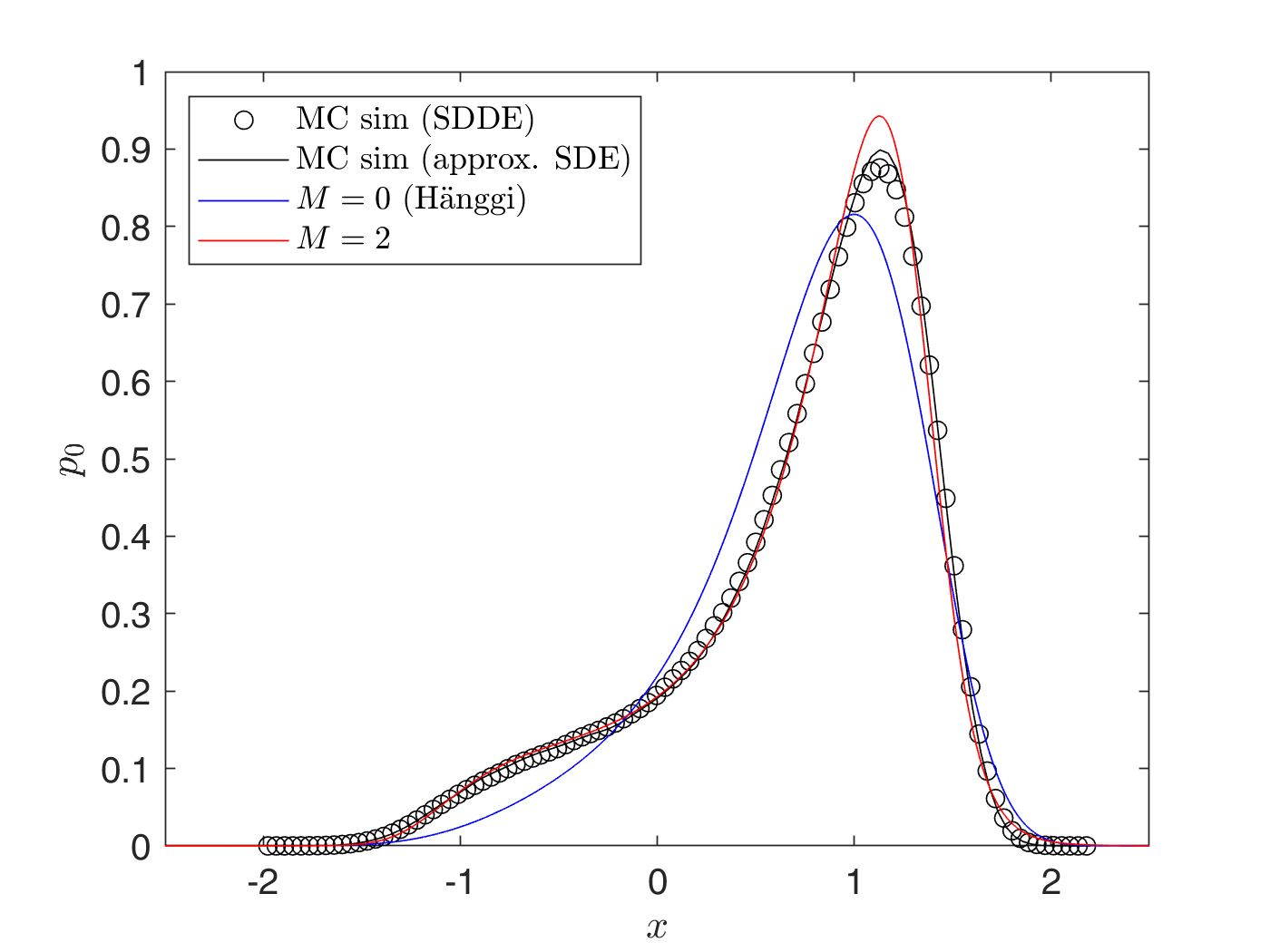}}
	\caption{Control of SDE~\eqref{eq:bistable_SDE} with $\sigma=1.2$, $s_{cor}=0.25$, $\hat{x}=1$, $a=1$, and increasing values of delay $\tau$. The uncontrolled bistable response is shown in (a). The controlled response is shown in (b) for $\tau=0.1$ ($\tilde{\sigma}=1.26$, $\tilde{s}_{cor}=0.28$), in (c) for $\tau=0.2$ ($\tilde{\sigma}=1.34$, $\tilde{s}_{cor}=0.31$), and in (d) for $\tau=0.4$ ($\tilde{\sigma}=1.55$, $\tilde{s}_{cor}=0.42$).} \label{fig:control1}
\end{figure*}

As a first example, we consider the SDE with the symmetric bistable potential
\begin{equation}
V(x)=\frac{x^4}{4}-\frac{x^2}{2}, 
\label{eq:bistable_potential}
\end{equation}
whose wells are located at $x=\pm1$. 
The resulting SDE, driven by additive OU excitation $\xi(t)$ reads
\begin{equation}
\frac{\id X(t)}{\id t}=-X^3(t)+X(t)+\sigma\xi(t). 
\label{eq:bistable_SDE}
\end{equation}
We designate the equilibrium located at $x_a=+1$ as the desirable equilibrium, and $x_b=-1$ as the undesirable one.  In this case, the controlled SDDE~\eqref{eq:SDDE} is given by
\begin{equation}
\frac{\id X(t)}{\id t}=-X^3(t)+X(t)-a(X(t-\tau)-\hat{x})+\sigma\xi(t). \label{eq:SDDE2}
\end{equation}
The effective potential $\tilde{V}(x)$ is given for this case by
\begin{equation}\label{eq:eff_potential1}
\tilde{V}(x)=V(x)+\frac{a}{2}(x-\hat{x})^2.
\end{equation}
By calculating its derivative, we see that the wells of effective potential are the roots of
\begin{equation}\label{eq:roots1}
\tilde{V}'(x)=x^3-(1-a)x-a\hat{x}=0.
\end{equation}
For the control gain $a=1$ and shift $\hat{x}=x_a=1$, the effective potential $\tilde{V}(x)$ has a single well at $x=1$. 
Figure~\ref{fig:control1} shows the corresponding stationary PDF for the control parameters $\{a=1,\hat x=1\}$ and several values of the delay time $\tau$.

The correlation time $s_{cor}$ of the OU noise and the time delay $\tau$ of the control are theoretically arbitrary. 
However, in applications, these time scales are often shorter than
the internal time scale of the system. We derive a characteristic time scale for system~\eqref{eq:bistable_SDE} based on its Lyapunov time. 
To this end, we linearize~\eqref{eq:bistable_SDE} for $\s=0$ around the desirable equilibrium $x_a=1$, to obtain the equation of variations, 
\begin{equation}\label{eq:lyapunov}
\frac{\id \delta X(t)}{\id t}=-V''(1)\delta X(t)=-2\delta X(t), \ \ \ \delta X(t_0)=\delta X_0,
\end{equation} 
for the perturbation $\delta X$, which admits the exact solution $\delta X(t)=\delta X_0e^{-2(t-t_0)}$. 
Therefore, the Lyapunov exponent is $\lambda=2$, leading to the characteristic Lyapunov time $\eta=\lambda^{-1}=0.5$. 
The time scale $\eta$ denotes the typical time it takes for a small perturbation to the equilibrium $x_a$ to decay. 
We allow the correlation time of the OU noise to be at most equal to the Lyapunov time scale $\eta$. 
More specifically, we present our results for several values of the correlation time in the interval $s_{cor}\in[0.1, 0.5]$.
Similarly, the control delay time $\tau$ is assumed to be strictly smaller than the characteristic time $\eta$.

However, we have to note that, while restricting $s_{cor}$ to be smaller that $\eta$ is a plausible physical argument, this does not imply that our response PDF approximation fails for larger values of correlation time. For the validity range of the approximation, see also the discussion in \cite[Sec. 4]{Mamis2019d}.

In the following subsections, we investigate various aspects of the controlled system and compare our 
approximate stationary PDF approach to the direct Monte Carlo simulations. 
Some of the discussed phenomena, such as peak drift, are unique to systems driven by colored noise
and are not present in the white noise case.

\subsubsection{Peak drift phenomenon}\label{sec:peak_drift}
In figure \ref{fig:control1}, the PDFs obtained by Monte Carlo simulations of both the uncontrolled SDE and the controlled SDDE, exhibit the so-called peak drift phenomenon, 
which had also been reported in, e.g., \cite[page 294]{Hanggi1995} and \cite[remark 4.1]{Mamis2019d}. 
Peak drift refers to the fact that the response PDF maxima are not observed at the wells of the potential $V(x)$, but they are slightly shifted. This phenomenon is only possible under colored noise excitation, and is not present in the case of additive white noise excitation. For white noise,  
we have $A_M=1/2$ in the PDF~\eqref{eq:p_0}, and therefore we can easily establish the equivalence,
\begin{equation}\label{equivalence}
\text{extrema of $\tilde{V}(x)$}\Leftrightarrow\text{extrema of $p_0(x)$}.
\end{equation}
Furthermore, since coefficient $A_0$, defined by Eq.~\eqref{eq:stationary_A_M}, is $x$-independent, we can easily see that equivalence~\eqref{equivalence} is also true for H{\"a}nggi's stationary PDF~\eqref{eq:p_0} for $M=0$. Thus, H{\"a}nggi's approximation of the stationary PDF $p_0$ fails to capture the observed peak drift phenomenon. This finding is also corroborated in figure \ref{fig:control1}. 

On the other hand, the critical points of our approximating stationary PDF~\eqref{eq:p_0} for $M=2$ are calculated as the roots of equation
\begin{equation}\label{eq:extrema2}
\tilde{V}'(x)+\tilde{\s}^2\tilde{A}'_2(x,R)=0,
\end{equation}
where $\tilde{A}_M(x,R)$ is defined by~\eqref{eq:stationary_A_M} for the effective quantities $\tilde{V}(x)$, $\tilde{\s}$, and $\tilde{s}_{cor}$.
Since $A_2'$ is generally non-zero, the peak of the stationary PDF does not necessarily coincide with the minima of the potential.
 Equation~\eqref{eq:extrema2} is the first analytic evidence of equivalence~\eqref{equivalence} not being true for additive colored excitations.

In the case of bistable potential~\eqref{eq:bistable_potential}, equation~\eqref{eq:extrema2} gives the depressed cubic equation,
\begin{equation}
c_3(R)x^3-c_1(R)x-a\hat{x}=0, 
\label{eq:extrema3}
\end{equation}
where
\begin{subequations}\label{eq:extrema_coeffs}
		\begin{equation}
	c_1(R)=(1-a)+\frac{3\tilde{\s}^2\tilde{s}_{cor}}{(1-\tilde{s}_{cor}R)^2}\left[1+\frac{2\tilde{s}_{cor}(1-a-R)}{1-\tilde{s}_{cor}R}\right],
	\end{equation}
	\begin{equation}
	c_3(R)=1+\frac{18\tilde{\s}^2\tilde{s}_{cor}^2}{(1-\tilde{s}_{cor}R)^3}.
	\end{equation} 
\end{subequations}
Note that Eq.~\eqref{eq:extrema3} also holds for the uncontrolled case, $a=0$, in which the tilded quantities are substituted by the untilded ones. Thus, a first straightforward consequence of~\eqref{eq:extrema3} is that, for the uncontrolled case, the extrema of the bimodal response  PDF are easily determined by the local minimum $x_0=0$, and the maxima 
\begin{equation}\label{eq:max_uncontrolled}
x_{1,2}=\pm\sqrt{\frac{1+\frac{3\s^2s_{cor}}{(1-s_{cor}R)^2}\left[1+\frac{2s_{cor}(1-R)}{1-s_{cor}R}\right]}{1+\frac{18\s^2s_{cor}^2}{(1-s_{cor}R)^3}}}.
\end{equation}
As expected, in the white noise limit, $s_{cor}=0$, Eq.~\eqref{eq:max_uncontrolled} results in $x_{1,2}=\pm1$, and so no peak drift is observed. 
Therefore, our approximating PDF, given by ~\eqref{eq:p_0} with $M=2$, allows for a systematic study of the peak drift phenomenon, without restoring to computationally expensive Monte Carlo simulations. 
\begin{figure*}
\subfigure[]{\label{fig:peak_a}\includegraphics[width=0.49\textwidth]{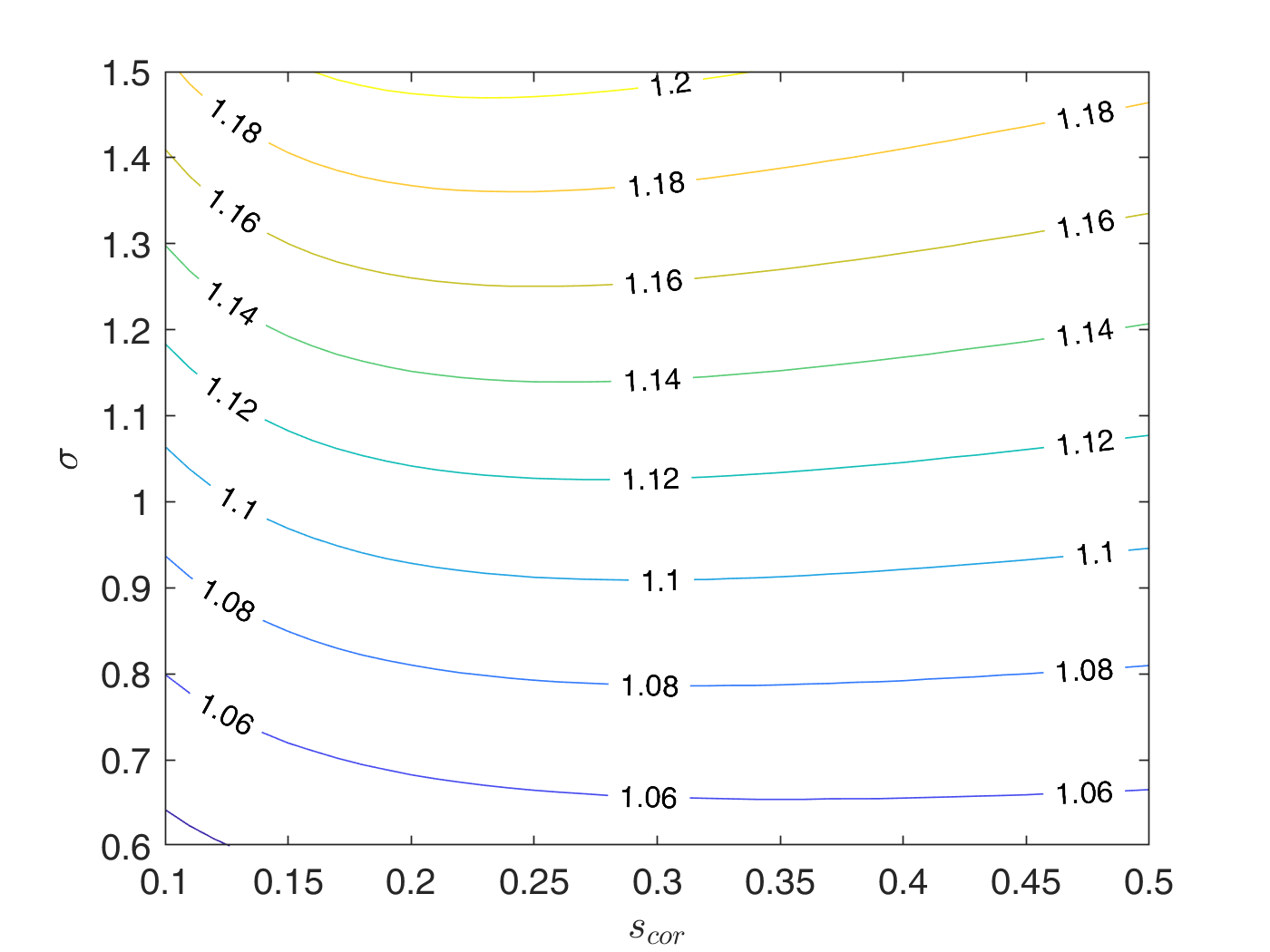}}
\subfigure[]{\label{fig:peak_b}\includegraphics[width=0.49\textwidth]{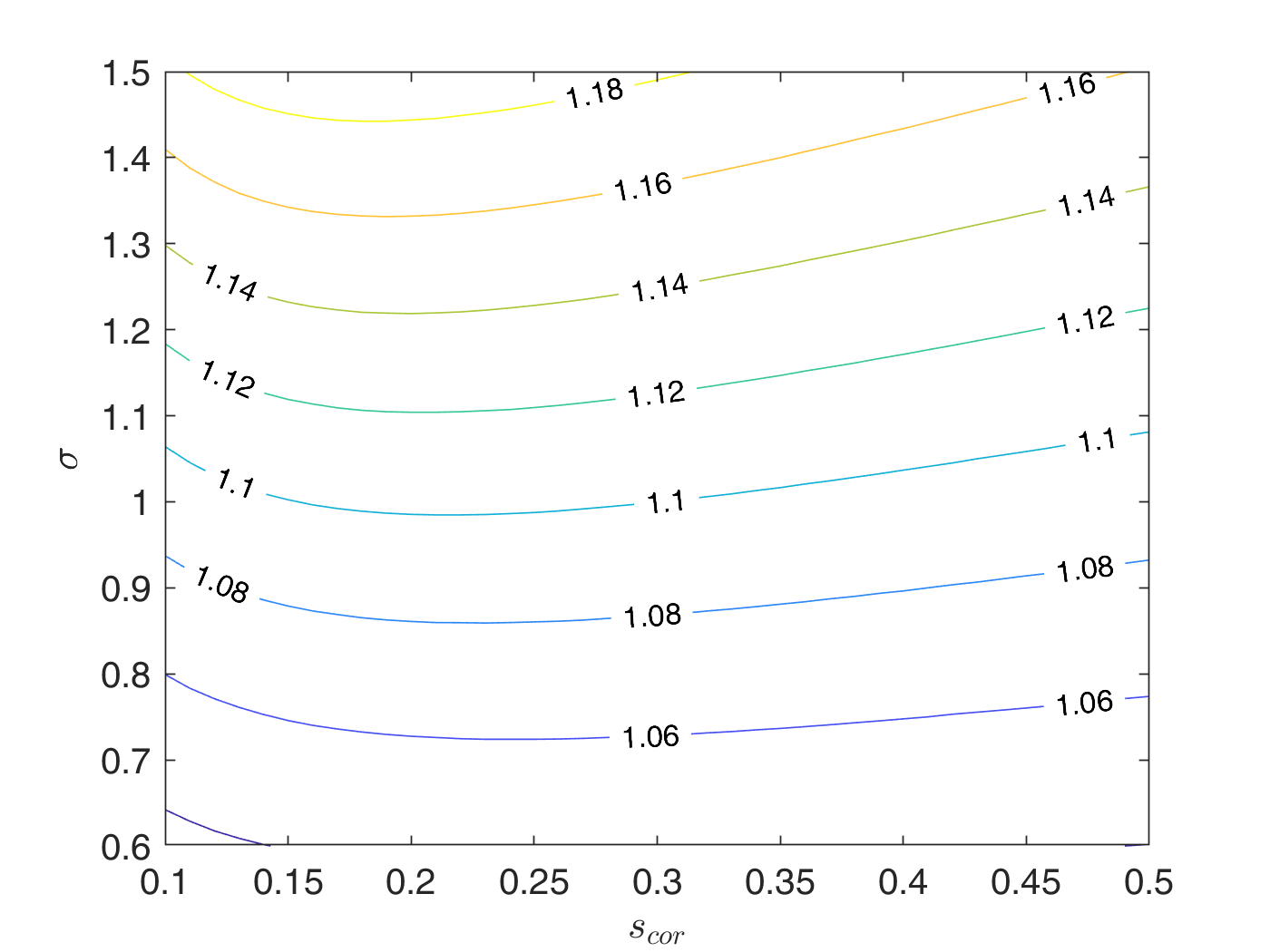}}
\subfigure[]{\label{fig:peak_c}\includegraphics[width=0.49\textwidth]{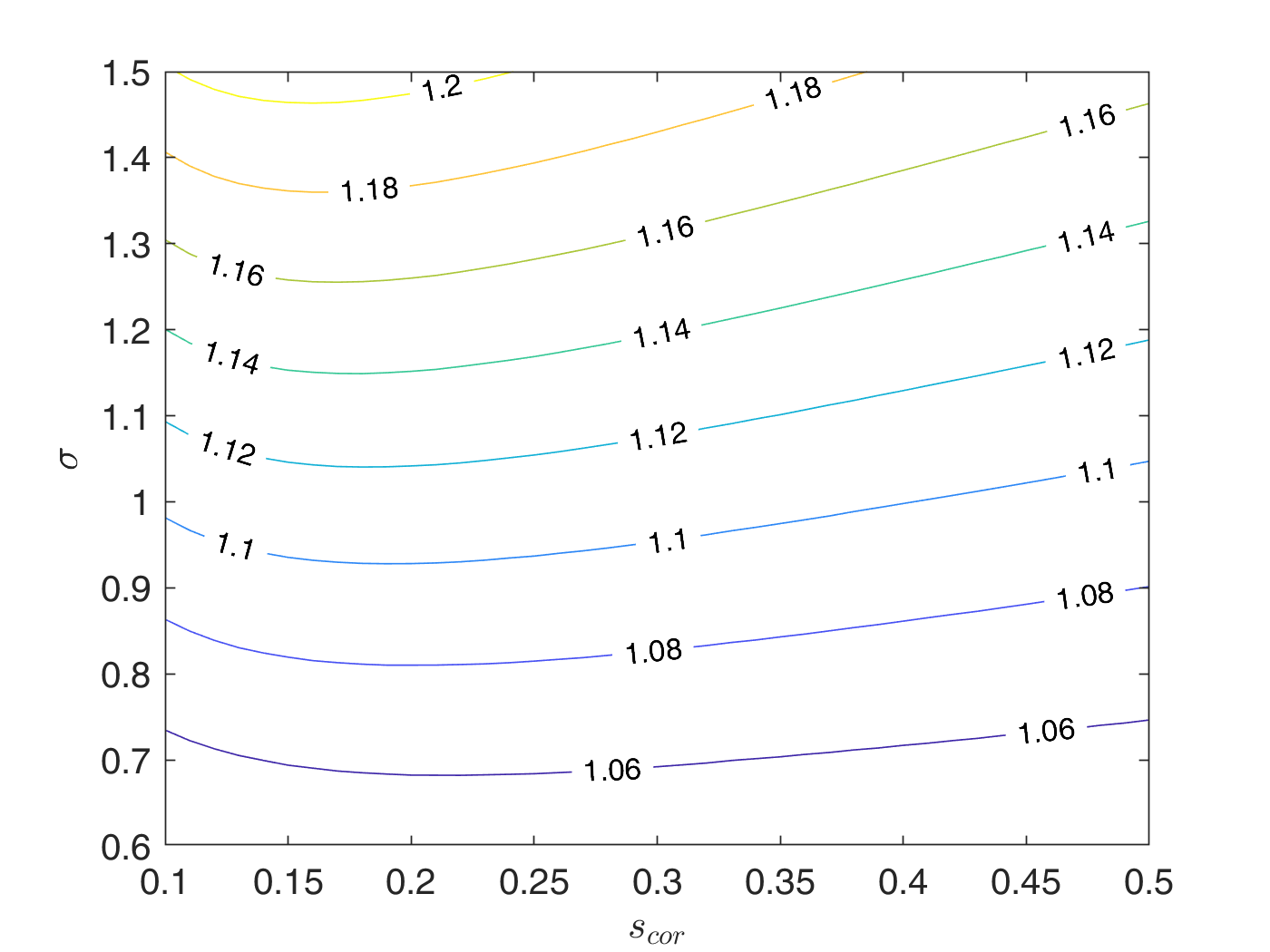}}
\subfigure[]{\label{fig:peak_d}\includegraphics[width=0.49\textwidth]{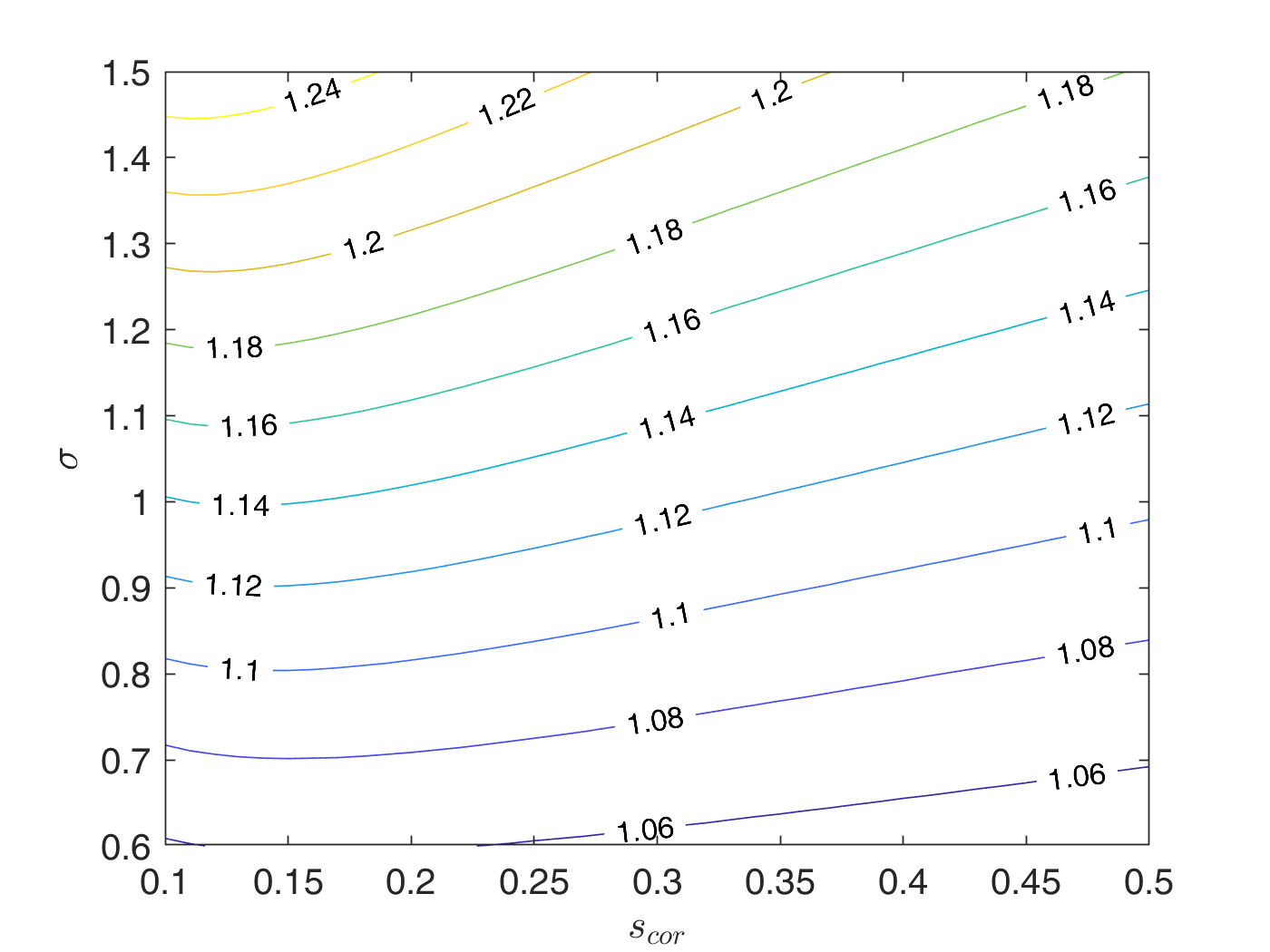}}
\caption{Contour plots of the peak $x$-coordinate for the uncontrolled and the controlled ($a=1$, $\hat{x}=1$) SDE~\eqref{eq:bistable_SDE}, as predicted by~\eqref{eq:p_0} for $M=2$. The uncontrolled case is shown in (a), while the controlled one is shown in (b) for $\tau=0.1$, in (c) for $\tau=0.2$, and in (d) for $\tau=0.4$. For the uncontrolled SDE, the response PDF is bimodal and symmetric around zero; therefore panel (a) shows the absolute value of the $x$-coordinates of the two peaks.} \label{fig:peak_drift}
\end{figure*}

Figure \ref{fig:peak_drift} shows the location of the peak of the stationary PDF for various parameter values. We first observe that the peak drift phenomenon is more pronounced in the bimodal response PDF of the uncontrolled SDE than the unimodal PDF of the controlled case. Also, increasing the noise intensity $\s$, as well as the time delay $\tau$ for the controlled case, increase the peak drift. On the other hand, the dependence of peak drift on the noise correlation time $s_{cor}$ is not monotone. For every pair of $\s$ and $\tau$, there is a value $s_{cor}$ for which the peak drift is maximized.  For larger values of $\s$ and $\tau$, the value $s_{cor}$ resulting in maximum peak drift is lowered. We note that the dependence of peak drift on the correlation time $s_{cor}$ is also briefly discussed in Ref.\cite[Sec. VI.B]{Hanggi1995}.

\subsubsection{Canceling the peak drift in the controlled system}\label{sec:canceling_peak_drift}
The above discussion shows that if the control shift $\hat x$ is chosen to be equal to the desirable equilibrium, i.e., $\hat x = x_a$, 
the resulting controlled PDF does not peak at $x_a$.
Therefore, the natural question is whether the control shift $\hat x$ can be chosen such that the peak drift is suppressed and the controlled response PDF attains its maximum at the desirable equilibrium $x_a$. The answer is yes.

Straightforward algebraic manipulations show that, in order for the depressed cubic equation~\eqref{eq:extrema3} to have $x_a$ as one of its roots, we must have
\begin{equation}\label{eq:x_hat}
\hat{x}(R)=\frac{x_a}{a}\left[c_3(R)x_a^2-c_1(R)\right],
\end{equation}
with coefficients $c_1(R)$ and $c_3(R)$ defined in~\eqref{eq:extrema_coeffs}. This control shift depends on the response moment $R$; and consequently, the effective potential also depend on $R$,
\begin{equation}\label{eq:eff_potential2}
\tilde{V}(x,R)=V(x)+\frac{a}{2}\left(x-\hat{x}(R)\right)^2.
\end{equation}
As a result, in the case of no peak drift control, the value $\hat{x}$ is a priori unknown, since it depends on the unknown response moment $R$. 
Nonetheless, we can substitute effective potential~\eqref{eq:eff_potential2} into the function $\mathcal{I}(R)$~\eqref{eq:I(R)} of the right-hand side of self-consistency Eq.~\eqref{eq:consistency}. Thus, for the present case, we specify $\mathcal{I}(R)$ into
\begin{equation}\label{eq:I(R)_2}
\mathcal{I}(R)=-\frac{\int_{\R}\frac{\tilde{V}''(x,R)}{\tilde{A}_M(x,R)}\exp\left(-\int^x\frac{1}{\tilde{\s}}\frac{\tilde{V}'(y,R)}{\tilde{A}_M(y,R)}\id y\right)\id x}{\int_{\R}\frac{1}{\tilde{A}_M(x,R)}\exp\left(-\int^x\frac{1}{\tilde{\s}}\frac{\tilde{V}'(y,R)}{\tilde{A}_M(y,R)}\id y\right)\id x}.
\end{equation}
By substituting $\mathcal{I}(R)$~\eqref{eq:I(R)_2} in self-consistency Eq.~\eqref{eq:consistency}, and following the iterative Algorithm~\ref{algorithm}, we can calculate the response moment $R$, and also the stationary response PDF. For the first step of Algorithm~\ref{algorithm}, the initial estimation of $R$ is calculated for $A_M=1/2$ and $\hat{x}=1$. Once the response moment is found, we are also able to calculate the value of appropriate control shift $\hat{x}(R)$ by~\eqref{eq:x_hat}. 

\begin{figure*}
\subfigure[]{\label{fig:peak_at_1_a}\includegraphics[width=0.49\textwidth]{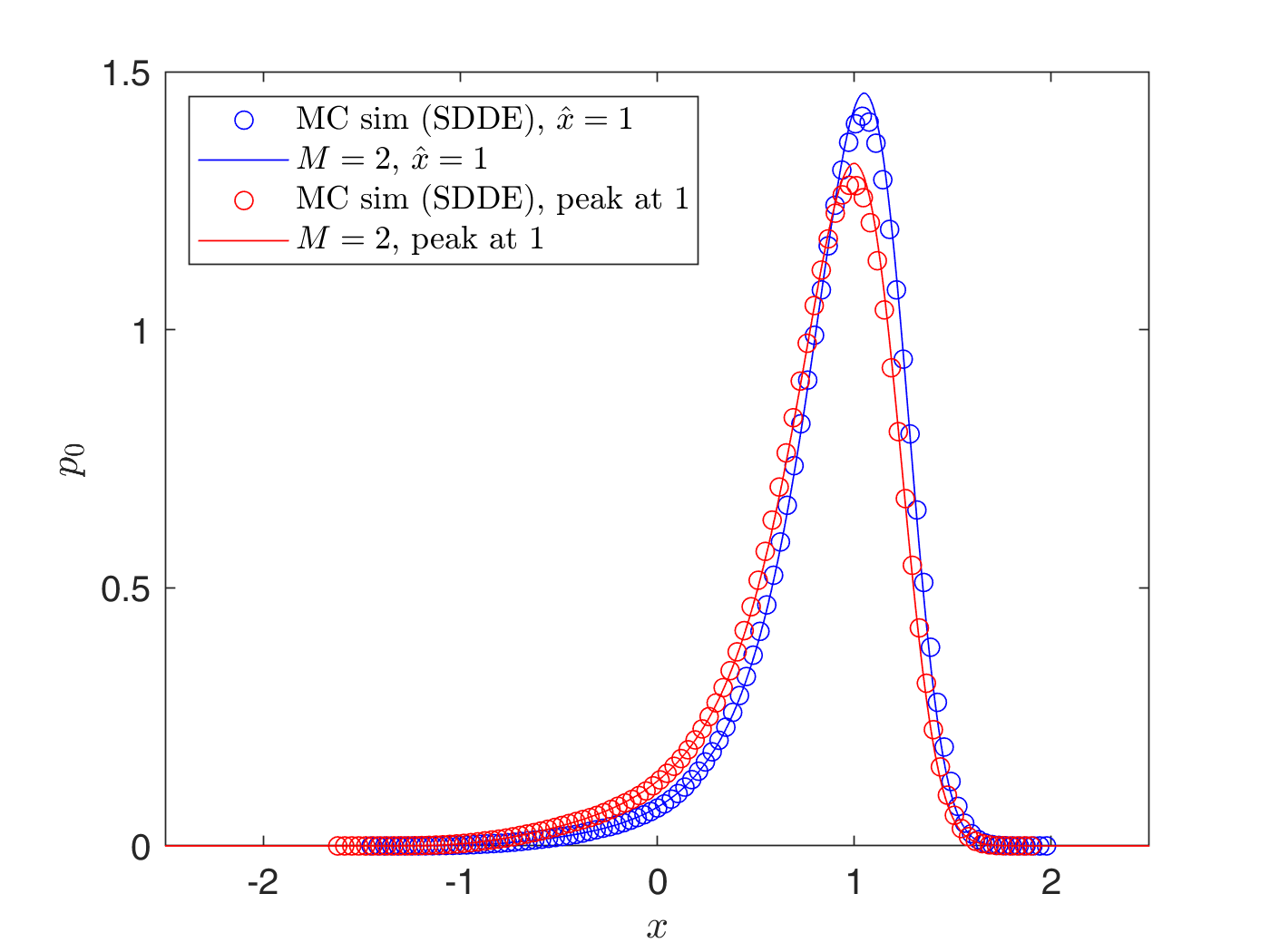}}
\subfigure[]{\label{fig:peak_at_1_b}\includegraphics[width=0.49\textwidth]{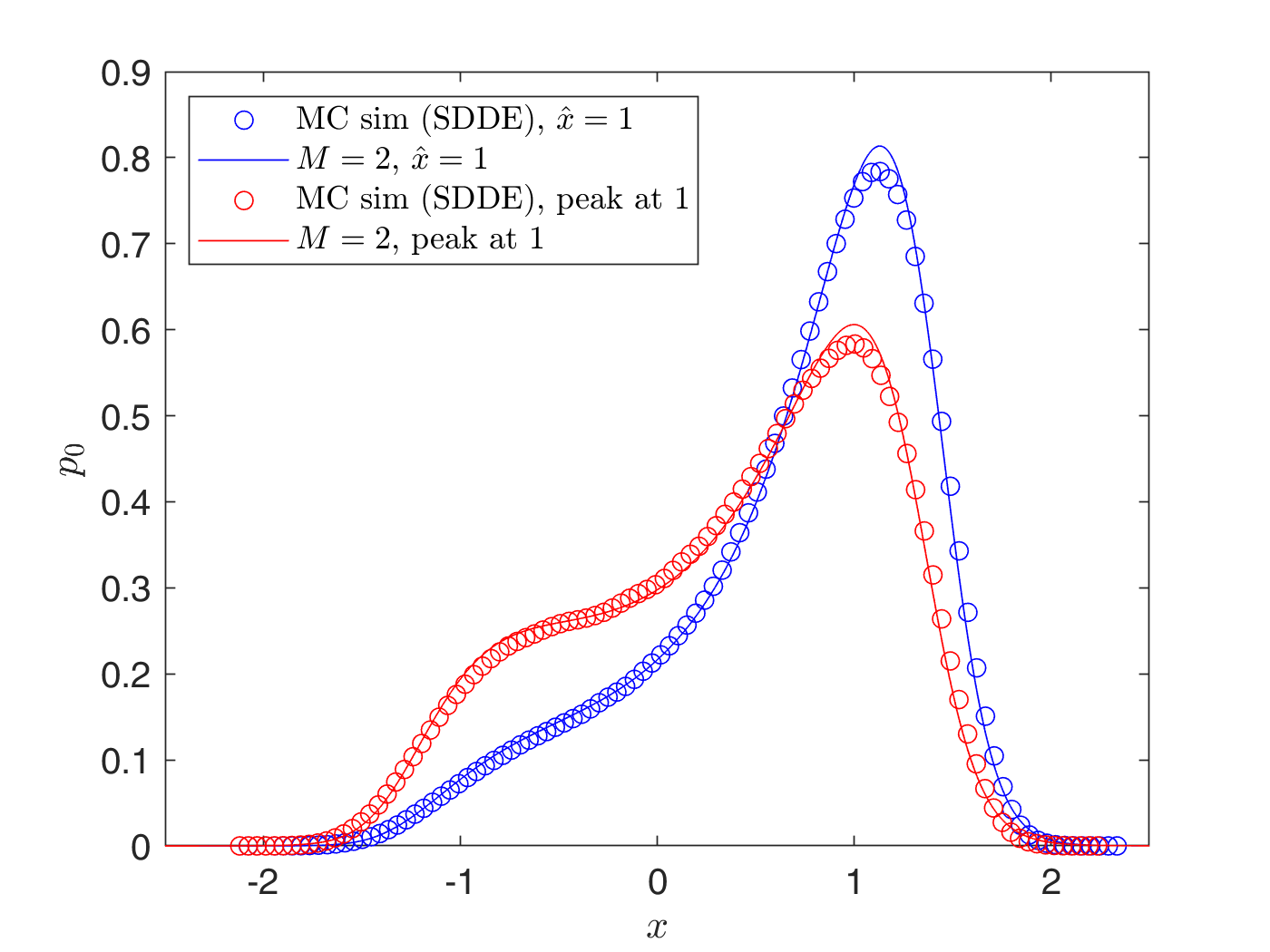}}
\caption{Control of SDE~\eqref{eq:bistable_SDE} with $a=1$, $s_{cor}=0.2$, $\tau=0.1$, $\s=0.8$ in (a), and $\s=1.4$ in (b), with $\hat{x}=1$ (blue curves), and with peak drift cancel (red curves). The $\hat{x}(R)$ in order to achieve peak drift canceling is calculated to 0.83 for (a), and 0.48 for (b). Stationary PDF approximation~\eqref{eq:p_0} for $M=2$ is plotted against the PDF obtained by Monte Carlo simulations of the respective SDDE, shown in plots by circles.} \label{fig:peak_at_1}
\end{figure*}

Figure \ref{fig:peak_at_1} shows the effectiveness of this procedure, resulting in controlled stationary PDFs with their peaks at 
the desirable equilibrium $x_a=1$. 
Note that the values of the control shifts $\hat{x}(R)$, required for canceling the peak drift, can be significantly smaller than $x_a=1$.

\subsubsection{Inflated tail in the controlled system}\label{sec:bump}
In figures \ref{fig:control1} and \ref{fig:peak_at_1_b}, we observe another feature of the controlled PDFs for systems excited by colored noise: the appearance of an inflation in the negative tail of the stationary PDF $p_0$ for certain combinations of noise $\s$, $s_{cor}$, and control parameters $a$, $\tau$. Our approximate PDF with $M=2$ accurately captures the inflated tail, whereas H{\"a}nggi's approximation completely misses these features. The formation of the inflated tail is due to the emergence of two additional inflection points in the negative tail of the unimodal controlled PDF. 
\begin{figure}
\centering
\includegraphics[width=0.49\textwidth]{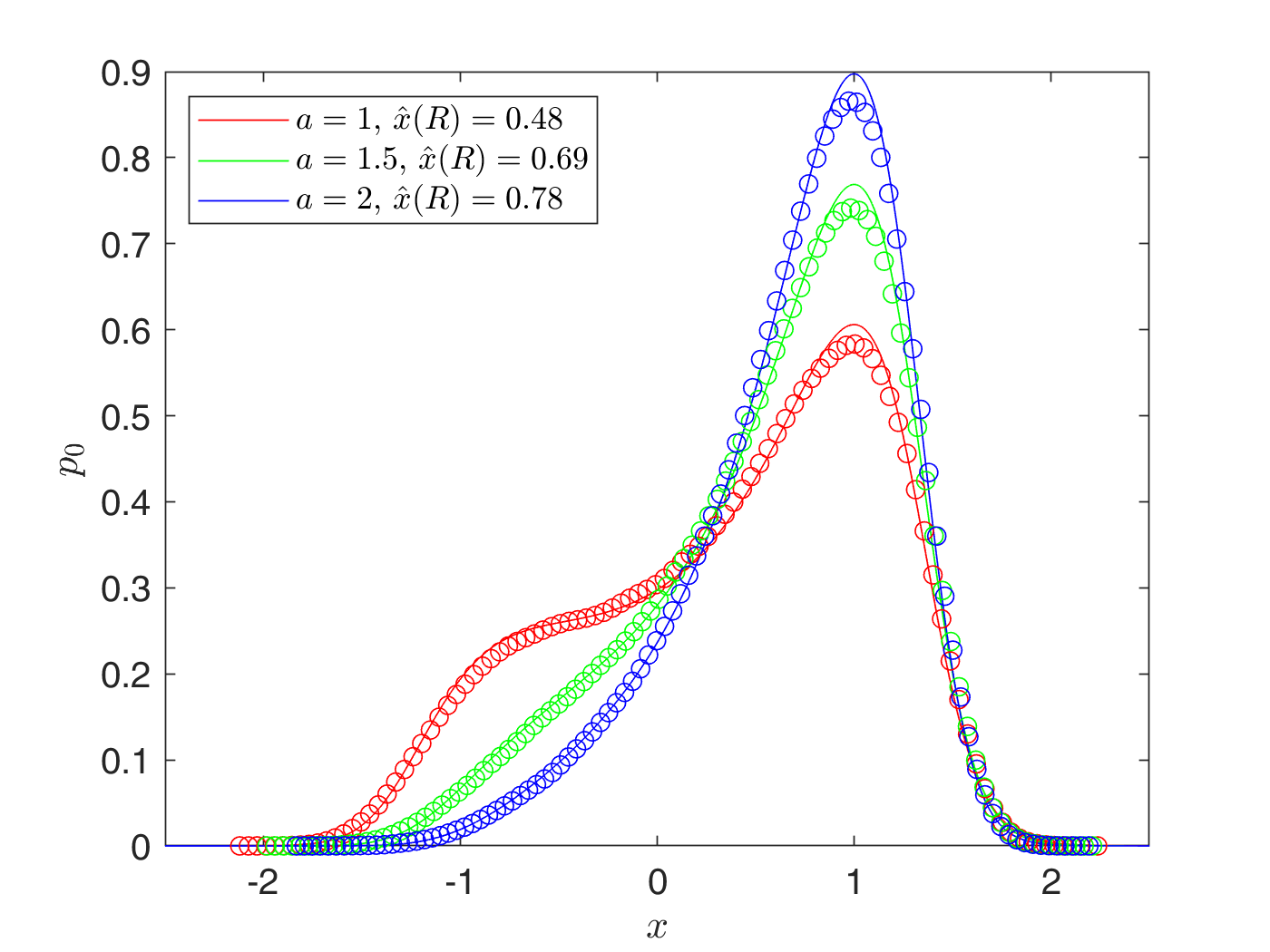}
\caption{Control of SDE~\eqref{eq:bistable_SDE} with $\s=1.4$, $s_{cor}=0.2$, $\tau=0.1$, no peak drift control, and increasing gain $a$.  Stationary PDF approximation~\eqref{eq:p_0} for $M=2$ (solid curves) is plotted against the PDF obtained by Monte Carlo simulations of SDDE~\eqref{eq:SDDE2}, shown by circles.}\label{fig:additive_control_gain}
\end{figure}

\begin{figure*}
\subfigure[]{\label{fig:pdf_form1}\includegraphics[width=0.49\textwidth]{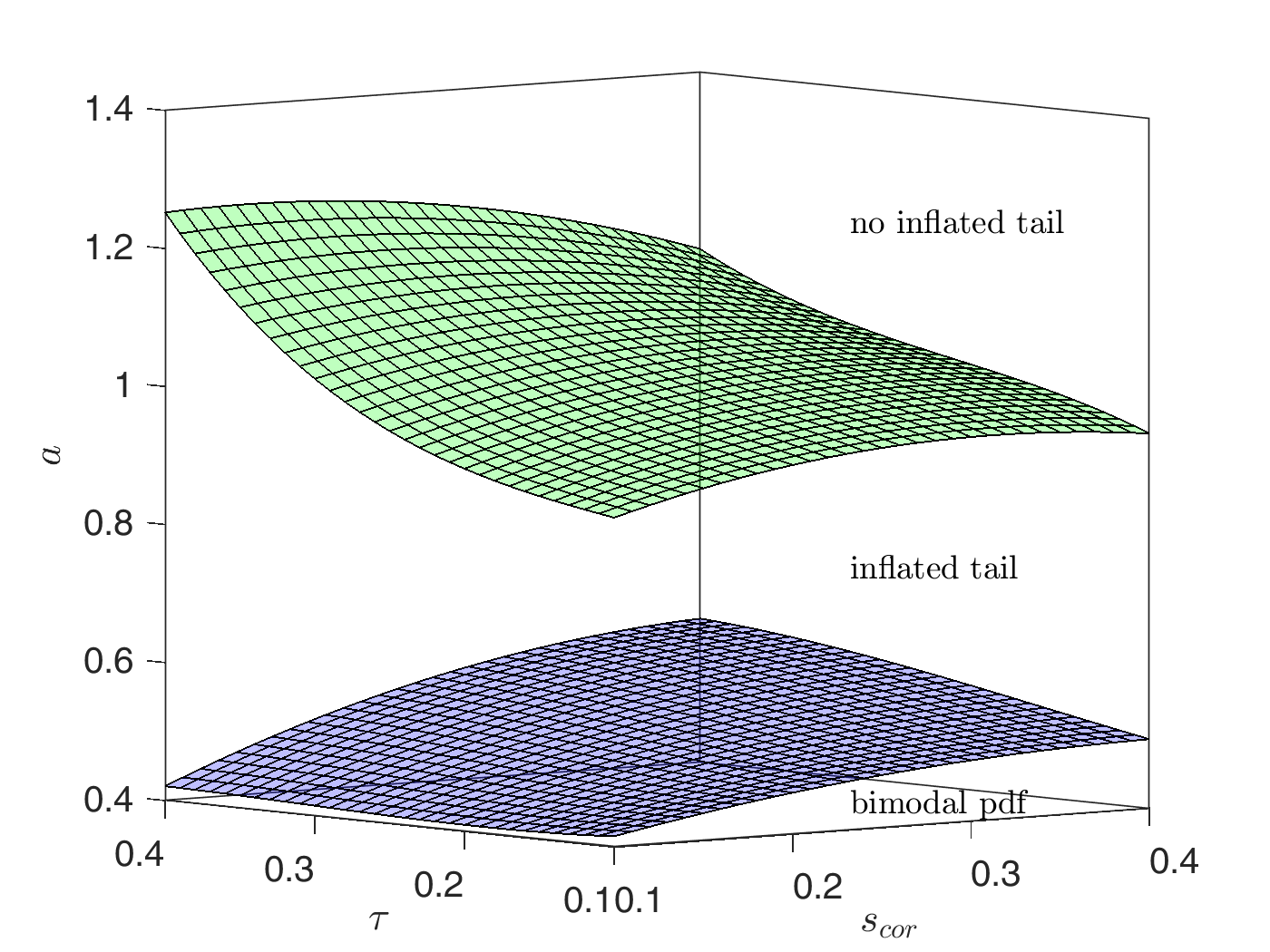}}
\subfigure[]{\label{fig:pdf_form2}\includegraphics[width=0.49\textwidth]{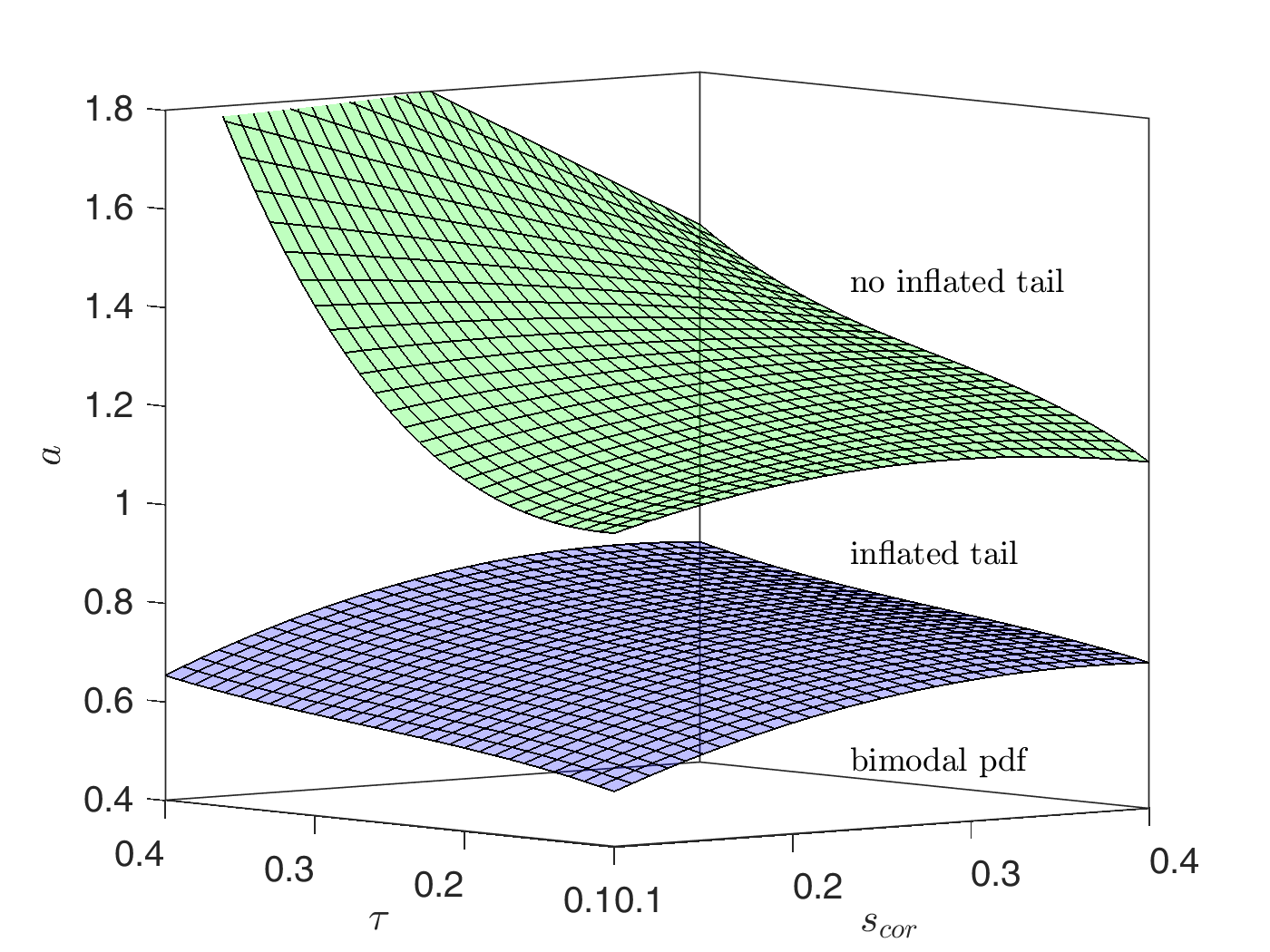}}
\caption{Three regimes of the controlled stationary PDF as a function of control parameters and correlation time of the noise.
The figure corresponds to SDE~\eqref{eq:bistable_SDE} with $\s=1$.
In (a) $\hat{x}=1$, and in (b) $\hat{x}$ is chosen such that the controlled PDF peak is located at $x=1$. In both figures, the region below the blue surface contains the parameter combinations that result in a bimodal PDF, the region between the two surfaces corresponds to unimodal PDFs with inflated tail, while the region above the green surface to unimodal PDFs with no inflated tail. } \label{fig:pdf_form}
\end{figure*}

A shown in figure~\ref{fig:additive_control_gain}, tail inflation can be suppressed by increasing the control gain $a$. 
Recall that the gain $a$ and the delay time $\tau$ must satisfy $a\tau<1$. As a result, the control gain cannot be arbitrarily large.
Therefore, there is a delicate balance between canceling the peak drift and suppressing rare transitions. 
In other words, given a control delay time $\tau$, the gain $a$ and the shift $\hat x$ should be chosen accordingly.
Since our approximate PDF approach quantifies the peak drift and the inflated tail, without requiring any computationally expensive simulations, the optimal control parameters can be readily identified by sweeping the parameter space.

Such parameter investigation results are shown in Fig.~\ref{fig:pdf_form}, for the case of $\s=1$. These results can be reproduced easily for any other value of noise intensity. In fig.~\ref{fig:pdf_form}, we see that, for fixed values of $s_{cor}$ and $\tau$ the response PDF is bimodal for small control gain $a$ values; this is the control failing to suppress transitions to the undesirable equilibrium. By increasing the value of $a$ we first obtain a unimodal PDF with inflated tail, and by further increase, a unimodal PDF with no tail inflation. Also, for larger values of time delay $\tau$, larger values of $a$ are required in order to obtain unimodal PDFs or suppression of the inflated tails. Last, as in peak drift phenomenon(see Fig.~\ref{fig:peak_drift}), the dependence of PDF form on $s_{cor}$ is not monotone. We also observe that, for the case of control with peak drift canceling, shown in Fig.~\ref{fig:pdf_form2}, the values of $a$ needed for unimodal PDF or for no tail inflation, are larger. Also, the dependence of separating surfaces on $\tau$ is more pronounced.  

\subsection{An optical laser excited by multiplicative noise} \label{sec:multiplicative}
A bistable SDE, arising in laser applications, is the model associated with the electromagnetically-induced transparency in a three-level atomic system inside an optical cavity.
Following Ref.~\cite{Wu2009}, the amplitude $X(t)$ of the transmitted intracavity light intensity is governed by
\begin{equation}\label{eq:SDE_multi}
\frac{\id X(t)}{\id t}=Y-c_1X(t)+c_2X^2(t)-c_3X^3(t)+\s X(t)\xi(t),
\end{equation}
where $Y$ is the incident light intensity amplitude. Reduced model~\eqref{eq:SDE_multi} is derived from the wave equation for the complex light intensity of the intracavity field, assuming rotating wave and slowly varying approximations. The two mechanisms in the wave equation that lead to the optical bistability are the nonlinear absorption, and the nonlinear refraction due to Kerr nonlinearity. 

\begin{figure}
\centering
\includegraphics[width=0.49\textwidth]{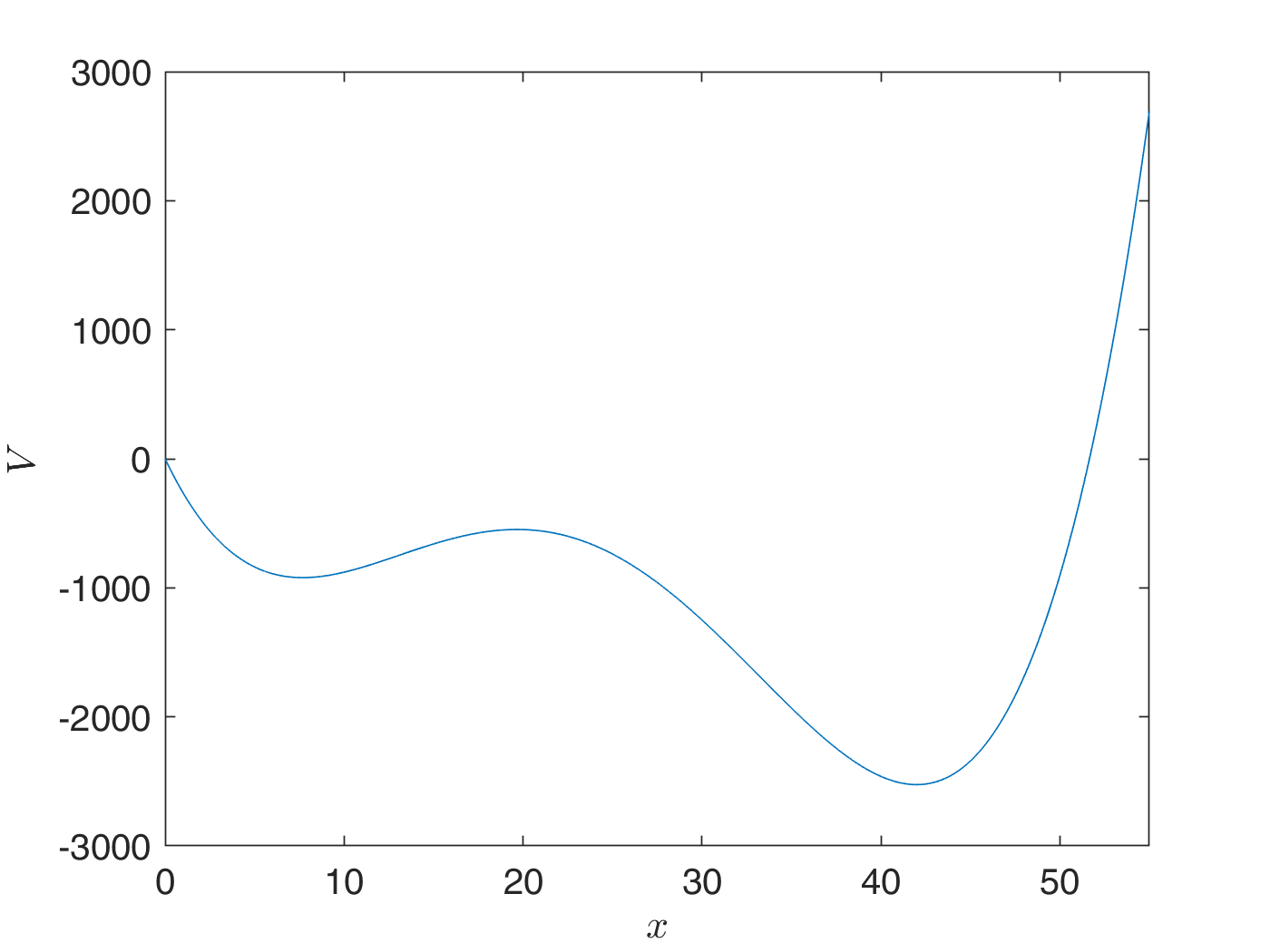}
\caption{Bistable potential~\eqref{eq:potential_multi} for the set of parameters~\eqref{eq:parameters}.}\label{fig:potential_multi}
\end{figure}

In this section, we use the same set of parameters considered in~\cite{Wu2009}, i.e.,
\begin{equation}\label{eq:parameters}
Y=292, \ c_1=59.79, \ c_2=3.19, \ c_3=0.046.
\end{equation}
The corresponding potential is given by
\begin{equation}\label{eq:potential_multi}
V(x)=\frac{c_3}{4}x^4-\frac{c_2}{3}x^3+\frac{c_1}{2}x^2-Yx,
\end{equation}
which has two minima at $x_a=42$ (global) and $x_b=7.69$ (local), and a local maximum at $x_0=19.66$ 
(see figure~\ref{fig:potential_multi}).
In SDE~\eqref{eq:SDE_multi}, the multiplicative excitation $\xi(t)$ models the noisy detuning of the cavity. Although in Ref.~\cite{Wu2009}, the white noise excitation is used for convenience, the authors mention that colored OU noise should be used in order to have a more realistic model. A similar emphasis on the relevance of colored noise is raised in the supplemental material of Ref.~\cite{Geng2020}. Therefore, here we consider the colored noise $\xi(t)$ obtained from the OU process.
Our numerical results are reported for the noise intensity $\s=2$, which is close to the noise intensities considered in \cite{Wu2009}. Furthermore, we use the correlation time $s_{cor}=0.02$, determined from the lower value of modulation frequency mentioned in the supplemental material of Ref.~\cite{Geng2020}.

\begin{figure*}
\subfigure[]{\label{fig:multi_a}\includegraphics[width=0.49\textwidth]{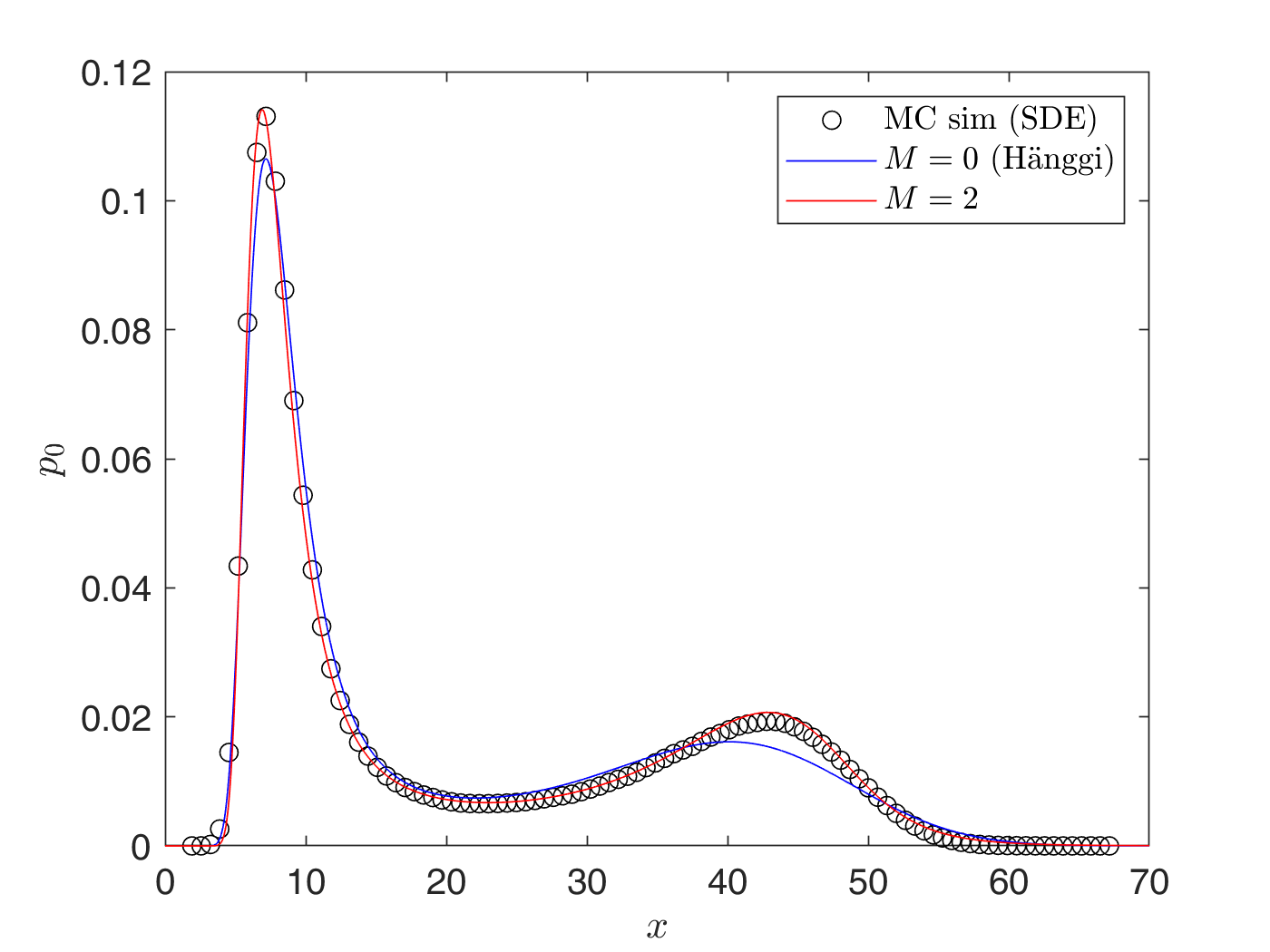}}
\subfigure[]{\label{fig:multi_b}\includegraphics[width=0.49\textwidth]{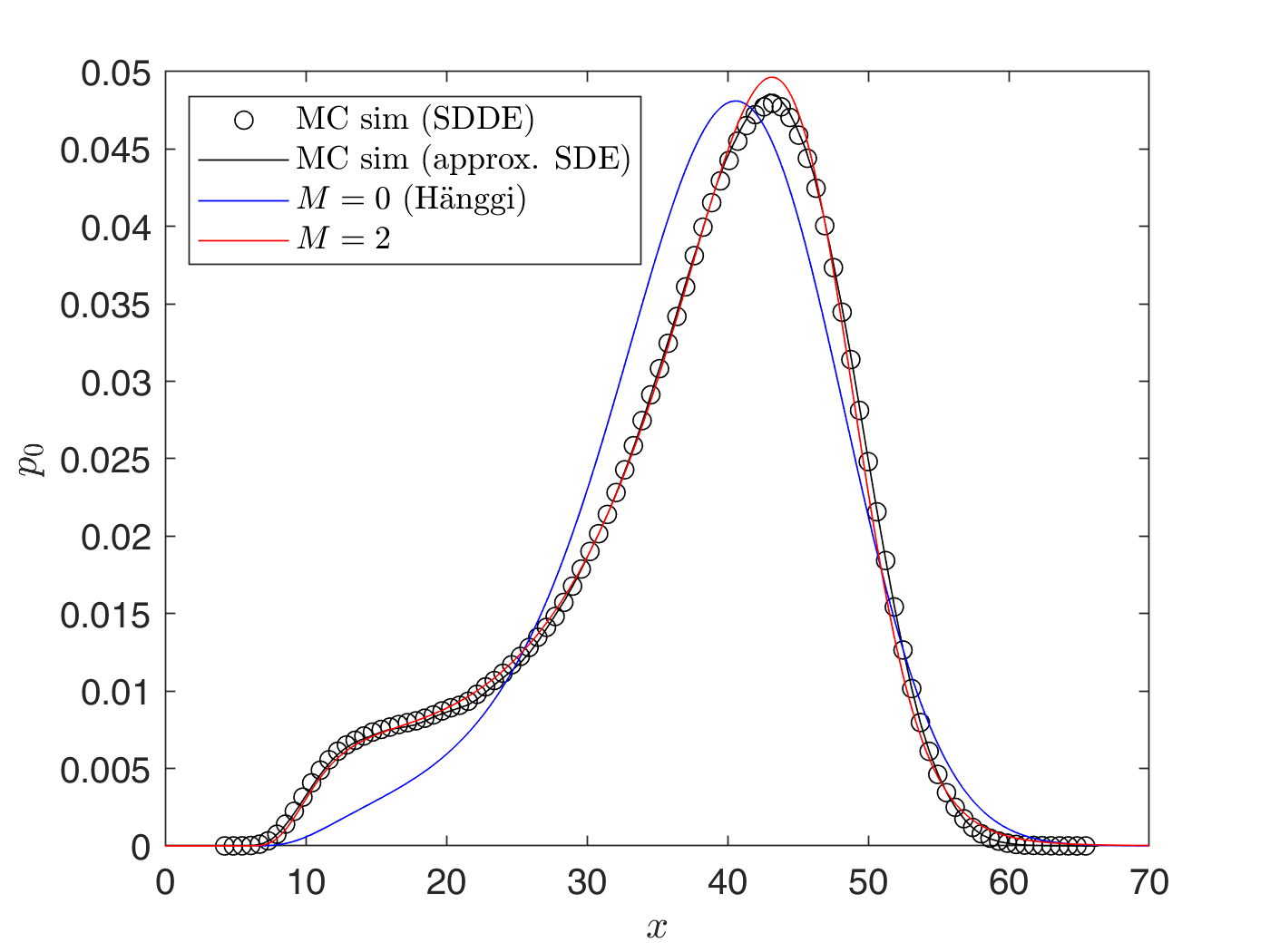}}
\subfigure[]{\label{fig:multi_c}\includegraphics[width=0.49\textwidth]{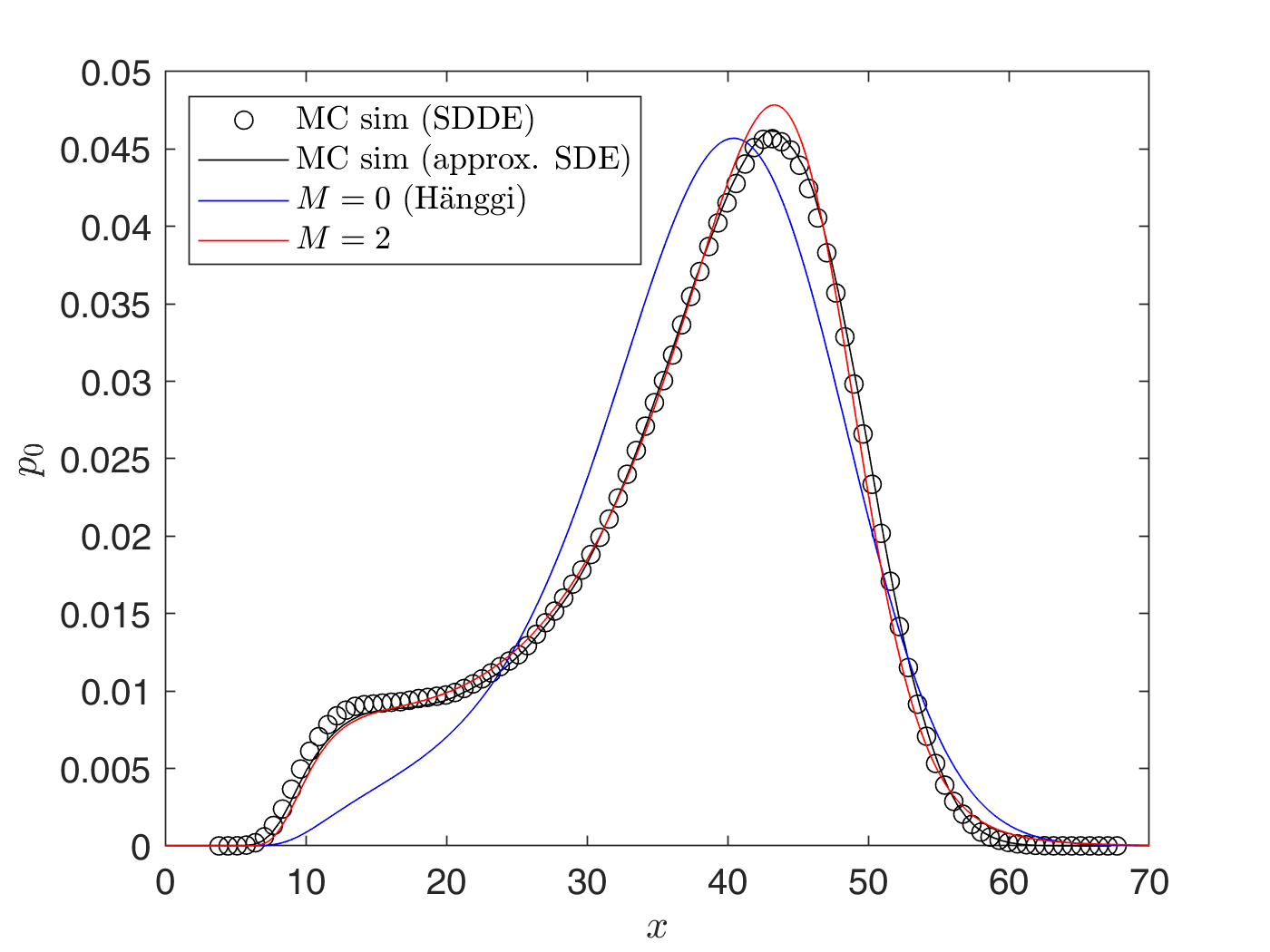}}
\subfigure[]{\label{fig:multi_d}\includegraphics[width=0.49\textwidth]{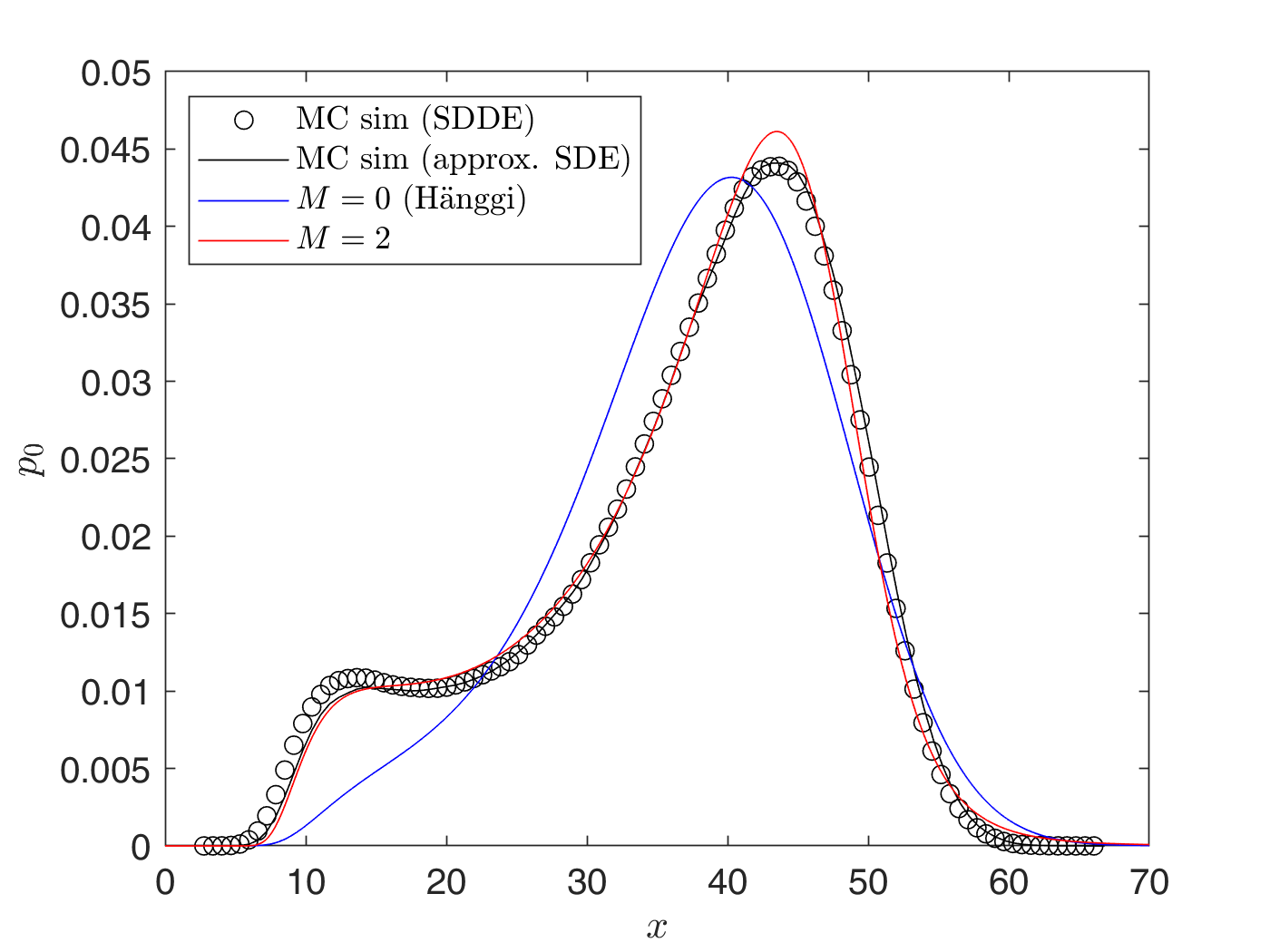}}
\caption{Control of SDE~\eqref{eq:SDE_multi} with $\hat{x}=42$, $a=4$, $\sigma=2$, $s_{cor}=0.02$, and increasing values of delay $\tau$. The uncontrolled bistable response is shown in (a). The controlled response is shown in (b) for $\tau=0.02$ ($\tilde{\sigma}=2.09$, $\tilde{s}_{cor}=0.022$), in (c) for $\tau=0.05$ ($\tilde{\sigma}=2.24$, $\tilde{s}_{cor}=0.025$), and in (d) for $\tau=0.08$ ($\tilde{\sigma}=2.43$, $\tilde{s}_{cor}=0.029$).} \label{fig:control_multi}
\end{figure*}

For the control of SDE~\eqref{eq:SDE_multi}, we choose the high-intensity equilibrium $x_a=42$ as the desirable one, so that the control SDDE reads
\begin{align}
\frac{\id X(t)}{\id t}=Y&-c_1X(t)+c_2X^2(t)-c_3X^3(t)+\nonumber\\&+a(X(t-\tau)-\hat{x})+\s X(t)\xi(t),\label{eq:SDDE_multi}
\end{align}
with the effective potential defined as
\begin{equation}\label{eq:eff_potential_multi}
\tilde{V}(x)=\frac{c_3}{4}x^4-\frac{c_2}{3}x^3+\frac{(c_1+a)}{2}x^2-(Y+a\hat{x})x+\frac{a}{2}\hat{x}^2.
\end{equation}
\begin{rem}[Peak drift phenomenon in the multiplicative case]
As in the additive case, the critical points of $p_0$ are easily determined as the roots of the equation,
\begin{equation}\label{eq:critical_multi}
\tilde{V}'(x)+\left(\tilde{\s}(x)\tilde{A}_M(x,R)\right)'\tilde{\s}(x)=0.
\end{equation}
This implies that, for a multiplicatively excited SDE, equivalence~\eqref{equivalence} between the wells of effective potential and response PDF maxima does not hold even in the white noise excitation case. More precisely, for the classical Fokker--Planck case, $A_M=1/2$, Eq.~\eqref{eq:critical_multi} reads
\begin{equation}\label{eq:roots_multi}
c_3x^3-c_2x^2+\left(c_1+a+\frac{\s^2}{2(1-a\tau)}\right)x-(Y+a\hat{x})=0.
\end{equation}
Thus, in the multiplicatively-excited case, the peak drift phenomenon is observed for both white and colored noise excitations.
\end{rem}
For $\s=2$ and $\hat{x}=x_a=42$, we choose $a=4$  as the control gain. For these values, Eq.~\eqref{eq:roots_multi} has a single root, and thus the controlled stationary PDF for white noise excitation is unimodal. Figure~\ref{fig:control_multi}
shows that the control effectively eliminates transitions to the undesirable equilibrium.

First, in figure \ref{fig:multi_a}, the bistable, uncontrolled response PDF is shown. We observe that its highest peak is around the local minimum ($x_b=7.69$) and not the global minimum ($x_a=42$) of the potential. This seemingly paradoxical observation is due to the multiplicative nature of the excitation; the $x$-dependent factor $|\s(x)|^{-1}$ in the right-hand side of PDF form~\eqref{eq:p_0} lowers the peak height around $x_a$.

In figures \ref{fig:multi_b}, \ref{fig:multi_c} and \ref{fig:multi_d}, several control delay times $\tau$ are chosen. In each case, our approximating PDF~\eqref{eq:p_0} with $M=2$ closely resembles the stationary response PDF obtained by direct Monte Carlo simulations. Furthermore, as in the additive case (see Section~\ref{sec:additive}), the inflated tail appears in the controlled PDFs, which our approximating PDF captures satisfactorily. On the other hand, H{\"a}nggi's approximation departs significantly from the Monte Carlo simulations and
completely misses the inflated tails.

\begin{figure}
\centering
\includegraphics[width=0.49\textwidth]{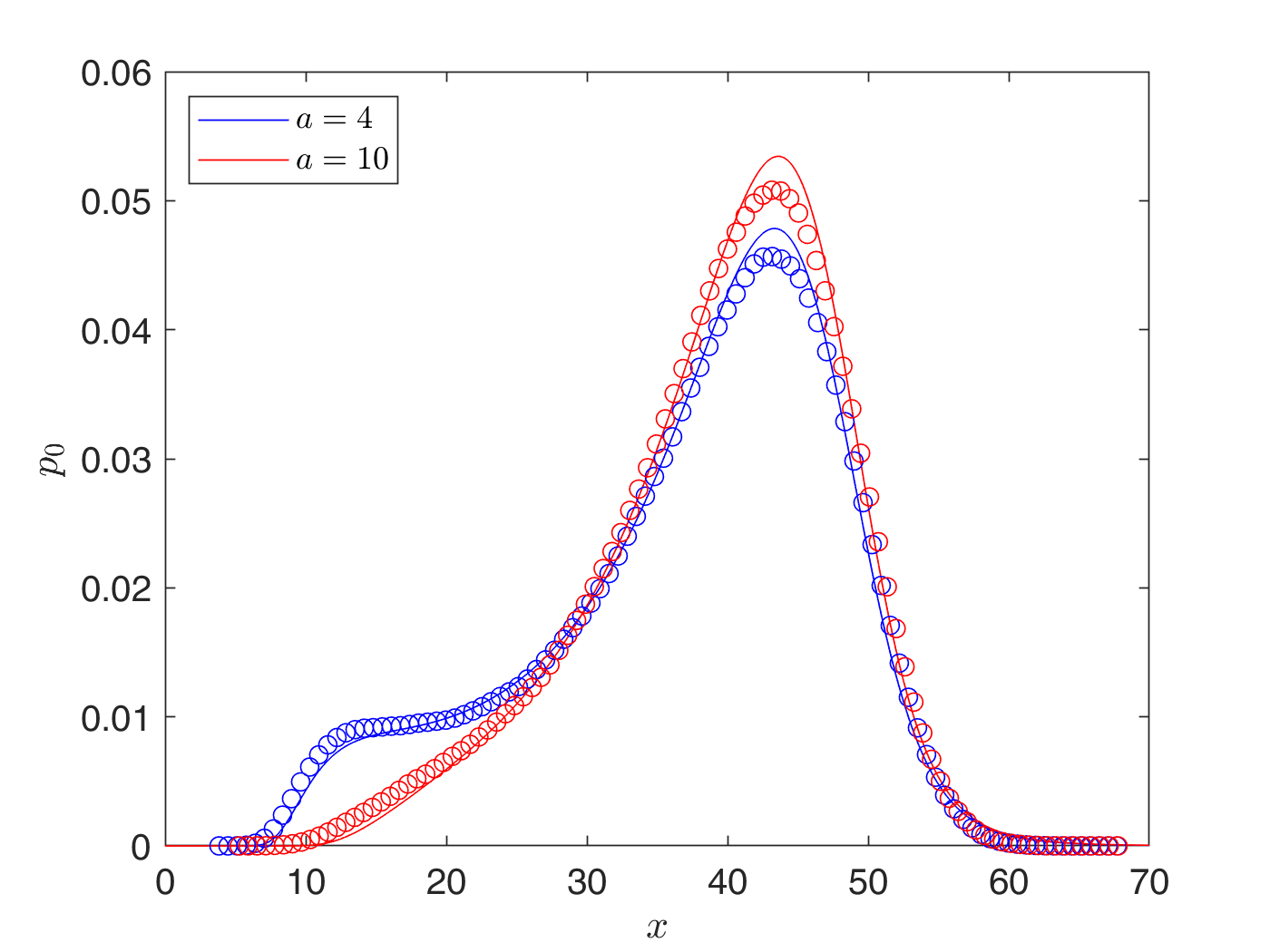}
\caption{Control of SDE~\eqref{eq:SDE_multi} with $\s=2$, $s_{cor}=0.02$, for $\tau=0.05$ and increasing control gain $a$.  Stationary PDF approximation~\eqref{eq:p_0} for $M=2$ (solid curves) is plotted against the PDF obtained by Monte Carlo simulations of SDDE~\eqref{eq:SDDE_multi}, shown by circles.}\label{fig:multi_control_gain}
\end{figure}

Finally, we observe that, for the relatively large time delay $\tau=0.08$ in figure~\ref{fig:multi_d}, a second undesirable peak begins to emerge around $x\simeq 10$. As shown in figure \ref{fig:multi_control_gain}, the tail inflation disappears by further increasing the control gain to $a=10$.

\section{Conclusions}\label{sec:conc}
We studied the mitigation of undesirable rare transitions in multistable SDEs excited by colored noise. The mitigation is  achieved by a time-delay feedback control that turns the original SDE into a SDDE. For small control delay $\tau$, the SDDE is approximated by an SDE. The approximating SDE reveals two competing effects of the controller:
(i) a stabilizing effect by deepening of the potential near the desirable equilibrium and (ii) a destabilizing effect by
effectively increasing the noise intensity.

In particular, for a scalar SDE excited by a colored OU process, the controller effectively 
increases both noise intensity $\s(x)$ and correlation time $s_{cor}$ by a factor of $1/\sqrt{1-a\tau}$ and $1/(1-a\tau)$, respectively. Here, $a$ is the control gain which cannot be arbitrary large since as $a$ approaches $1/\tau$, 
the effective noise intensity grows indefinitely. As a result, the control parameters are not arbitrary and need to be chosen 
judiciously.

We proposed a parsimonious method for choosing the optimal control parameters that guarantee the mitigation of 
undesirable rare transitions. Our method relies on a  nonlinear Fokker--Planck equation that was recently derived by Mamis et al.~\cite{Mamis2019d} for SDEs with additive noise. Here, we generalized this equation to 
the case of multiplicative noise. This nonlinear Fokker--Planck equation governs the evolution of the response PDF of a stochastic equation, and its stationary solution can be estimated by a rapidly convergent iterative algorithm.
As a result, we are able to sweep the control parameter space and select the optimal control parameters, 
without requiring any expensive Monte Carlo simulations.

We demonstrated the efficacy of our method on two examples: a bistable additively excited SDE, which serves as a classical benchmark problem~\cite{Hanggi1995, Ridolfi2011}, and a multiplicatively excited SDE modeling optical bistabilities in 
lasers~\cite{Wu2009, Peters2019}.  In each case, our approximate PDFs are in excellent agreement with the true PDF of the controlled SDDE constructed by direct Monte Carlo simulations. In particular, our method captures two features of the PDF that the widely-used  H{\"a}nggi's approximation~\cite{Hanggi1995} fails to predict. These are the peak drift and the inflated tail of the stationary PDFs, as discussion in sections~\ref{sec:peak_drift} and~\ref{sec:bump}. 

There are two important directions for future work. The first one is controlling multidimensional stochastic dynamical systems excited by colored noise. This is in principle achievable since  the multidimensional analog of the scalar nonlinear Fokker--Planck equations used here can be formulated, see e.g. \cite{Mamis2018, Mamis2020}. The main difficulty in this direction is the fact that analytical stationary solutions for multidimensional nonlinear Fokker--Planck equations are often unavailable, even for relatively simple systems such as stochastic oscillators~\cite{Soize1994, Zhu1992, Zhu2001, Mamis2016e}. Thus, in the multidimensional case, an efficient numerical solver has to be employed to approximate the solution of the nonlinear Fokker--Planck equation. Whether obtaining this numerical solution is less computationally expensive than the direct Monte Carlo simulations will be the main question regarding the efficacy of our method in the multidimensional case. In any case however, the formulation and validation of accurate nonlinear Fokker--Planck equations will be a useful semi-analytic tool for the study of multidimensional systems under colored noise excitation.

Another interesting direction is to derive the appropriate nonlinear Fokker--Planck equations that correspond to the controlled SDDE, without the small time-delay limitation. The derivation of such equations is laborious yet feasible, with the preparatory work concerning the relevant functional analysis questions having already been performed~\cite{Mamis2020}.

\appendix

\section{Proof of Theorem \ref{th_rescaling}} \label{A_rescaling}
Assuming that $a\tau<1$, consider $t$ as a function of rescaled time $s$, so that $t(s)=(1-a\tau)s$. First we define the stochastic processes $\bm{X}(s):=\bm{X}(t(s))$ and $\bm{\xi}(s):=\bm{\xi}(t(s))$. Since $\id t=(1-a\tau)\id s$, approximating SDE~\eqref{eq:Ndim_approx_SDE} is expressed equivalently as
\begin{equation}\label{eq:A1}
\frac{\id \bm{X}(s)}{\id s}=\frac{\id \bm{X}(t)}{\id t}(1-a\tau)=-\nabla \tilde{V}(\bm{X}(s))+\bm{\s}(\bm{X}(s))\bm{\xi}(s),
\end{equation}
which is the desired rescaled approximating SDE~\eqref{eq:Ndim_resc_SDE}. The SDE for rescaled noise $\bm{\xi}(s)$ is derived from filter SDE~\eqref{eq:filter} as
\begin{align}
\frac{\id \bm{\xi}(s)}{\id s}=\frac{\id \bm{\xi}(t)}{\id t}(1-a\tau)&=(1-a\tau)\bm{\alpha}(\bm{\xi}(s))+\nonumber\\&+(1-a\tau)\bm{\beta}(\bm{\xi}(s))\bm{\xi}^{\text{WN}}(t).\label{eq:A3}
\end{align}
In Eq.~\eqref{eq:A3}, the white noise $\bm{\xi}^{\text{WN}}(t)$ should also be transformed to the rescaled time $s$. Thus,  we introduce $\tilde{\bm{\xi}}^{\text{WN}}(s)=\sqrt{1-a\tau}\bm{\xi}^{\text{WN}}(t)$ and show that it is a standard Gaussian white noise. By its definition, $\tilde{\bm{\xi}}^{\text{WN}}(s)$ is a Gaussian, zero-mean process, and its autocorrelation is calculated to 
\begin{align*}
&\mathbb{E}\left[\tilde{\bm{\xi}}^{\text{WN}}(s_1)\left(\tilde{\bm{\xi}}^{\text{WN}}(s_2)\right)^T\right]=\nonumber\\ &=(1-a\tau)\mathbb{E}\left[\bm{\xi}^{\text{WN}}(t_1)\left(\bm{\xi}^{\text{WN}}(t_2)\right)^T\right]=\bm{I}(1-a\tau)\delta(t_1-t_2)\nonumber\\ &=\bm{I}(1-a\tau)\delta\left((1-a\tau)(s_1-s_2)\right)=\bm{I}\delta(s_1-s_2),
\end{align*} 
with the last step employing the delta function scale property. This proves that $\tilde{\bm{\xi}}^{\text{WN}}(s)$ is a standard white noise, and its substitution in rescaled filter SDE~\eqref{eq:A3} reads
\begin{equation}\label{eq:A4}
\frac{\id \bm{\xi}(s)}{\id s}=(1-a\tau)\bm{\alpha}(\bm{\xi}(s))+\sqrt{1-a\tau}\bm{\beta}(\bm{\xi}(s))\tilde{\bm{\xi}}^{\text{WN}}(s).
\end{equation}
By omitting the tilde from rescaled white noise, we obtain the desired rescaled filter SDE~\eqref{eq:Ndim_resc_filter}.

\section{Proof of Corollary \ref{cor:rescaled_OU}}\label{A:cor}
By applying theorem \ref{th_rescaling}, the system consisting of approximating SDE~\eqref{eq:Ndim_approx_SDE} and SDE filter producing the standard OU noise with autocorrelation~\eqref{eq:OU_multi_cor2}, is rescaled to 
\begin{equation}\label{eq:cor1}
\begin{aligned}
&\frac{\id \bm{X}(s)}{\id s}=-\nabla\tilde{V}(\bm{X}(s))+\bm{\s}(\bm{X}_(s))\bm{X}(s)\\
&\frac{\id \bm{\xi}(s)}{\id s}=-(1-a\tau)\bm{A}\bm{\xi}(s)+\sqrt{1-a\tau}\bm{A}\bm{\xi}^{\text{WN}}(s),
\end{aligned}
\end{equation}
where $\bm{A}=\text{diag}\left[\left\{1/s_{cor}^{(\ell)}\right\}_{\ell=1}^m\right]$, see Eq.~\eqref{eq:A_B}. By introducing $\tilde{\bm{\xi}}(s)=\sqrt{1-a\tau}\bm{\xi}(s)$, rescaled SDE~\eqref{eq:cor1} is expressed equivalently as
\begin{subequations}
\begin{align}
&\frac{\id \bm{X}(s)}{\id s}=-\nabla\tilde{V}(\bm{X}(s))+\frac{\bm{\s}(\bm{X}(s))}{\sqrt{1-a\tau}}\tilde{\xi}(s)\label{eq:cor2a}\\
&\frac{\id \tilde{\bm{\xi}}(s)}{\id s}=-(1-a\tau)\bm{A}\tilde{\bm{\xi}}(s)+(1-a\tau)\bm{A}\bm{\xi}^{\text{WN}}(s).\label{eq:cor2b}
\end{align}
\end{subequations}
By SDE~\eqref{eq:cor2b}, $\tilde{\bm{\xi}}(s)$ is identified as a standard OU noise with drift matrix $\tilde{\bm{A}}=(1-a\tau)\bm{A}=\text{diag}\left[\left\{(1-a\tau)/s_{cor}^{(\ell)}\right\}_{\ell=1}^m\right]$. By introducing the tilded noise intensity~\eqref{eq:tilde_sigma} and correlation times~\eqref{eq:tilde_scor}, we obtain the desired rescaled approximating SDE~\eqref{eq:Ndim_SDE_OU_aug} under standard OU noise.

\section{Derivation of nonlinear Fokker--Planck equations} \label{A:1}
The derivation of nonlinear Fokker--Planck equations begins with the stochastic Liouville equation (also called the colored noise master equation) corresponding to stochastic dynamical systems under colored noise. In \cite[Sec. III.D]{Hanggi1995}, stochastic Liouville equation for SDE~\eqref{eq:SDE} is formulated to  
\begin{widetext}
\begin{equation}
\frac{\partial p(x,t)}{\partial t}-\frac{\partial}{\partial x}\left[V'(x)p(x,t)\right]=\frac{\partial}{\partial x}\s(x)\frac{\partial}{\partial x}\s(x)\int_{t_0}^tC_{\xi}(t,s)\mathbb{E}\left[\delta(X(t)-x)\exp\int_s^t\zeta(X(u))\id u\right]\id s, \label{eq:SLE}
\end{equation}
\end{widetext}
with $\zeta(x)$ defined by Eq.~\eqref{eq:zeta1}. For the derivation of SLE~\eqref{eq:SLE}, we employ what H{\"a}nggi in \cite{Hanggi1989} calls the Fox's trick \cite{Fox1986a}. Fox's trick is essentially the substitution of noise excitation $\xi(t)$ via SDE~\eqref{eq:SDE}. In order for this substitution to be legitimate, noise intensity $\s(x)$ has to be non-vanishing, hence the assumption $\s(x)\neq 0$ at the last paragraph of Sec. \ref{sec:controlled_dynamical_systems}.

Stochastic Liouville Eq.~\eqref{eq:SLE} is an exact evolution equation for PDF $p(x,t)$, but it is not closed, due to the presence of the average in its right-hand side. In it, the random delta function $\delta(X(t)-x)$ appears, whose defining property reads
\begin{equation}
\mathbb{E}\left[\delta(X(t)-x)\right]=\int_{\R}\delta(y-x)p(y,t)\id y=p(x,t).
\end{equation}
Following \cite{Mamis2019d}, and in order to obtain a nonlinear Fokker--Planck equation in closed form from Eq.~\eqref{eq:SLE}, we apply a current-time approximation to the exponential in the right-hand side of~\eqref{eq:SLE}. First, we decompose the integrand into its mean value $\mathbb{E}\left[\zeta(X(u))\right]$ and the fluctuation around it $\phi(x,u)=\zeta(x)-\mathbb{E}\left[\zeta(X(u))\right]$:
\begin{align} 
&\exp\int_s^t\zeta(X(u))\id u=\nonumber\\&=\exp\int_s^t\mathbb{E}\left[\zeta(X(u))\right]\id u\cdot\exp\int_s^t\phi(X(u),u)\id u.\label{eq:decomposition}
\end{align}
Then, by assuming that $\phi(X(u),u)$ is small, we approximate the fluctuations exponential by a quadratic Taylor expansion with respect to $s$ around current time $t$. The first and second temporal derivatives of fluctuations exponential are easily calculated to
\begin{align*}
\frac{\partial}{\partial s}\exp&\int_s^t\phi(X(u),u)\id u=\nonumber\\&=-\phi(X(s),s)\exp\int_s^t\phi(X(u),u)\id u,
\end{align*}
\begin{align*}
&\frac{\partial^2}{\partial s^2}\exp\int_s^t\phi(X(u),u)\id u=\nonumber\\&=\left[\phi^2(X(s),s)-\frac{\partial\phi(X(s),s)}{\partial s}\right]\exp\int_s^t\phi(X(u),u)\id u.
\end{align*}
By further assuming that the first temporal derivative of the fluctuations is also small, we obtain the following current-time approximation for the whole exponential term in Eq.~\eqref{eq:SLE}:
\begin{align}
&\exp\int_s^t\zeta(X(u))\id u\approx\nonumber\\&\approx\exp\int_s^t\mathbb{E}\left[\zeta(X(u))\right]\id u\cdot\sum_{m=0}^2\frac{\phi^m(X(t),t)(t-s)^m}{m!}. \label{eq:approximation}
\end{align}
By substituting the above current-time approximation into the stochastic Liouville equation, we obtain the nonlinear Fokker--Planck equation:
\begin{equation} \label{eq:pdf_equation1}
\frac{\partial p(x,t)}{\partial t}-\frac{\partial}{\partial x}\left[V'(x)p(x,t)\right]=\frac{\partial}{\partial x}\s(x)\frac{\partial}{\partial x}\s(x) A_2(x,t) p(x,t),
\end{equation}
with
\begin{equation}
A_2(x,t)=\sum_{m=0}^2\frac{D_m(t)}{m!}\left\{\zeta(x)-\mathbb{E}\left[\zeta(X(t))\right]\right\}^m
\end{equation}
and
\begin{equation}
D_m(t)=\int_{t_0}^tC_{\xi}(t,s)\exp\int_s^t\mathbb{E}\left[\zeta(X(u))\right]\id u\ (t-s)^m\id s.
\end{equation}
By easy algebraic manipulations on the right-hand side of~\eqref{eq:pdf_equation1}, Eq.~\eqref{eq:nFokker--Planck} for $M=2$ is obtained. If, in the approximation scheme~\eqref{eq:approximation}, only the zeroth order term in Taylor series is employed, the usual H{\"a}nggi's Eq.~\eqref{eq:nFokker--Planck} for $M=0$ is retrieved. 

%\bibliographystyle{unsrtnat}
%\bibliography{../library,../bibliog}

\begin{thebibliography}{65}
	\providecommand{\natexlab}[1]{#1}
	\providecommand{\url}[1]{\texttt{#1}}
	\expandafter\ifx\csname urlstyle\endcsname\relax
	\providecommand{\doi}[1]{doi: #1}\else
	\providecommand{\doi}{doi: \begingroup \urlstyle{rm}\Url}\fi
	
	\bibitem[Dakos et~al.(2008)Dakos, Scheffer, {Van Nes}, Brovkin, Petoukhov, and
	Held]{Dakos2008}
	V.~Dakos, M.~Scheffer, E.~H. {Van Nes}, V.~Brovkin, V.~Petoukhov, and H.~Held.
	\newblock {Slowing down as an early warning signal for abrupt climate change}.
	\newblock \emph{Proceedings of the National Academy of Sciences of the United
		States of America}, 105\penalty0 (38):\penalty0 14308--14312, sep 2008.
	\newblock ISSN 00278424.
	\newblock \doi{10.1073/pnas.0802430105}.
	\newblock URL \url{www.pnas.org/cgi/content/full/}.
	
	\bibitem[Lenton et~al.(2008)Lenton, Held, Kriegler, Hall, Lucht, Rahmstorf, and
	Schellnhuber]{Lenton2008}
	T.~M. Lenton, H.~Held, E.~Kriegler, J.~W. Hall, W.~Lucht, S.~Rahmstorf, and
	H.~J. Schellnhuber.
	\newblock {Tipping elements in the Earth's climate system}, feb 2008.
	\newblock ISSN 00278424.
	\newblock URL \url{www.pnas.org/cgi/content/full/}.
	
	\bibitem[Mendez and Farazmand(2020)]{Mendez2020}
	A.~Mendez and M.~Farazmand.
	\newblock {Mitigating climate tipping points under various emission reduction
		and carbon capture scenarios (preprint)}.
	\newblock \emph{arXiv:2012.01613}, dec 2020.
	\newblock URL \url{http://arxiv.org/abs/2012.01613}.
	
	\bibitem[Zhu and Zhu(2010)]{Zhu2010}
	P.~Zhu and Y.~J. Zhu.
	\newblock {Statistical Properties of Intensity Fluctuation of Saturation Laser
		Model Driven By Cross-Correlated Additive and Multiplicative Noises}.
	\newblock \emph{International Journal of Modern Physics B}, 24\penalty0
	(14):\penalty0 2175--2188, 2010.
	\newblock ISSN 0217-9792.
	\newblock \doi{10.1142/S0217979210055755}.
	\newblock URL
	\url{http://www.worldscientific.com/doi/abs/10.1142/S0217979210055755}.
	
	\bibitem[Ridolfi et~al.(2011)Ridolfi, D'Odorico, and Lalo]{Ridolfi2011}
	L.~Ridolfi, P.~D'Odorico, and F.~Lalo.
	\newblock \emph{{Noise-Induced Phenomena in the Environmental Sciences}}.
	\newblock Cambridge University Press, 2011.
	
	\bibitem[Spanio et~al.(2017)Spanio, Hidalgo, and Mu{\~{n}}oz]{Spanio2017}
	T.~Spanio, J.~Hidalgo, and M.~A. Mu{\~{n}}oz.
	\newblock {Impact of environmental colored noise in single-species population
		dynamics}.
	\newblock \emph{Physical Review E}, 96\penalty0 (4):\penalty0 1--9, 2017.
	\newblock ISSN 24700053.
	\newblock \doi{10.1103/PhysRevE.96.042301}.
	
	\bibitem[Zeng et~al.(2017)Zeng, Xie, Wang, Zhang, Dong, Guan, Li, and
	Duan]{Zeng2017b}
	C.~Zeng, Q.~Xie, T.~Wang, C.~Zhang, X.~Dong, L.~Guan, K.~Li, and W.~Duan.
	\newblock {Stochastic ecological kinetics of regime shifts in a time-delayed
		lake eutrophication ecosystem}.
	\newblock \emph{Ecosphere}, 8\penalty0 (6), 2017.
	\newblock ISSN 21508925.
	\newblock \doi{10.1002/ecs2.1805}.
	
	\bibitem[Bose and Trimper(2011)]{Bose2011}
	T.~Bose and S.~Trimper.
	\newblock {Noise-assisted interactions of tumor and immune cells}.
	\newblock \emph{Physical Review E}, 84\penalty0 (021927):\penalty0 1--9, 2011.
	
	\bibitem[Idris and {Abu Bakar}(2016)]{Idris2016}
	I.~M. Idris and M.~R. {Abu Bakar}.
	\newblock {Effect of tumor microenvironmental factors on tumor growth dynamics
		modeled by correlated colored noises with colored cross-correlation}.
	\newblock \emph{Physica A: Statistical Mechanics and its Applications},
	453:\penalty0 298--304, 2016.
	\newblock ISSN 03784371.
	\newblock \doi{10.1016/j.physa.2016.01.082}.
	\newblock URL \url{http://dx.doi.org/10.1016/j.physa.2016.01.082}.
	
	\bibitem[Yang et~al.(2014)Yang, Han, Zeng, Wang, Liu, Zhang, and
	Tian]{Yang2014}
	T.~Yang, Q.~L. Han, C.~H. Zeng, H.~Wang, Z.~Q. Liu, C.~Zhang, and D.~Tian.
	\newblock {Transition and resonance induced by colored noises in tumor model
		under immune surveillance}.
	\newblock \emph{Indian Journal of Physics}, 88\penalty0 (11):\penalty0
	1211--1219, 2014.
	\newblock ISSN 09749845.
	\newblock \doi{10.1007/s12648-014-0521-7}.
	
	\bibitem[Zeng and Wang(2010)]{Zeng2010}
	C.~Zeng and H.~Wang.
	\newblock {Colored Noise Enhanced Stability in a Tumor Cell Growth System Under
		Immune Response}.
	\newblock \emph{Journal of Statistical Physics}, 141\penalty0 (5):\penalty0
	889--908, 2010.
	\newblock ISSN 00224715.
	\newblock \doi{10.1007/s10955-010-0068-8}.
	
	\bibitem[Li and Ning(2016)]{Li2016}
	X.~L. Li and L.~J. Ning.
	\newblock {Effect of correlation in FitzHugh–Nagumo model with non-Gaussian
		noise and multiplicative signal}.
	\newblock \emph{Indian Journal of Physics}, 90\penalty0 (1):\penalty0 91--98,
	2016.
	\newblock ISSN 0973-1458.
	\newblock \doi{10.1007/s12648-015-0717-5}.
	\newblock URL \url{http://link.springer.com/10.1007/s12648-015-0717-5}.
	
	\bibitem[Li and Zhu(2018)]{Li2018}
	S.~H. Li and Q.~X. Zhu.
	\newblock {Stochastic impact in Fitzhugh–Nagumo neural system with time
		delays driven by colored noises}.
	\newblock \emph{Chinese Journal of Physics}, 56\penalty0 (1):\penalty0
	346--354, 2018.
	\newblock ISSN 05779073.
	\newblock \doi{10.1016/j.cjph.2017.11.014}.
	
	\bibitem[Bose and Trimper(2012)]{Bose2012}
	T.~Bose and S.~Trimper.
	\newblock {Influence of randomness and retardation on the FMR-linewidth}.
	\newblock \emph{Physica Status Solidi (B) Basic Research}, 249\penalty0
	(1):\penalty0 172--180, 2012.
	\newblock ISSN 03701972.
	\newblock \doi{10.1002/pssb.201147164}.
	
	\bibitem[Chattopadhyay and Aifantis(2016)]{Chattopadhyay2016}
	A.~K. Chattopadhyay and E.~C. Aifantis.
	\newblock {Stochastically forced dislocation density distribution in plastic
		deformation}.
	\newblock \emph{Physical Review E}, 94\penalty0 (2):\penalty0 1--7, 2016.
	\newblock ISSN 24700053.
	\newblock \doi{10.1103/PhysRevE.94.022139}.
	
	\bibitem[Gayout et~al.(2021)Gayout, Bourgoin, and Plihon]{Gayout2021}
	A.~Gayout, M.~Bourgoin, and N.~Plihon.
	\newblock {Rare Event-Triggered Transitions in Aerodynamic Bifurcation}.
	\newblock \emph{Physical Review Letters}, 126\penalty0 (10):\penalty0 104501,
	mar 2021.
	\newblock ISSN 0031-9007.
	\newblock \doi{10.1103/physrevlett.126.104501}.
	\newblock URL
	\url{https://journals.aps.org/prl/abstract/10.1103/PhysRevLett.126.104501}.
	
	\bibitem[Dallas et~al.(2020)Dallas, Seshasayanan, and Fauve]{Dallas2020}
	V.~Dallas, K.~Seshasayanan, and S.~Fauve.
	\newblock {Transitions between turbulent states in a two-dimensional shear
		flow}.
	\newblock \emph{Physical Review Fluids}, 5\penalty0 (8):\penalty0 084610, aug
	2020.
	\newblock ISSN 2469990X.
	\newblock \doi{10.1103/PhysRevFluids.5.084610}.
	\newblock URL
	\url{https://journals.aps.org/prfluids/abstract/10.1103/PhysRevFluids.5.084610}.
	
	\bibitem[{Van Kan} et~al.(2019){Van Kan}, Nemoto, and Alexakis]{VanKan2019}
	A.~{Van Kan}, T.~Nemoto, and A.~Alexakis.
	\newblock {Rare transitions to thin-layer turbulent condensates}.
	\newblock \emph{Journal of Fluid Mechanics}, 878:\penalty0 356--369, nov 2019.
	\newblock ISSN 14697645.
	\newblock \doi{10.1017/jfm.2019.572}.
	\newblock URL
	\url{https://www.cambridge.org/core/journals/journal-of-fluid-mechanics/article/abs/rare-transitions-to-thinlayer-turbulent-condensates/5F30B3ED0485F779F10A28B1E6D539EC}.
	
	\bibitem[Shukla et~al.(2016)Shukla, Fauve, and Brachet]{Shukla2016}
	V.~Shukla, S.~Fauve, and M.~Brachet.
	\newblock {Statistical theory of reversals in two-dimensional confined
		turbulent flows}.
	\newblock \emph{Physical Review E}, 94\penalty0 (6):\penalty0 061101, dec 2016.
	\newblock ISSN 24700053.
	\newblock \doi{10.1103/PhysRevE.94.061101}.
	\newblock URL
	\url{https://journals.aps.org/pre/abstract/10.1103/PhysRevE.94.061101}.
	
	\bibitem[Farazmand(2016)]{faraz_adjoint}
	M.~Farazmand.
	\newblock An adjoint-based approach for finding invariant solutions of
	{N}avier-{S}tokes equations.
	\newblock \emph{J. Fluid Mech.}, 795:\penalty0 278--312, 2016.
	\newblock \doi{10.1017/jfm.2016.203}.
	
	\bibitem[Farazmand and Sapsis(2017)]{Farazmand2017}
	M.~Farazmand and T.~P. Sapsis.
	\newblock {A variational approach to probing extreme events in turbulent
		dynamical systems}.
	\newblock \emph{Science Advances}, 3\penalty0 (9):\penalty0 e1701533, sep 2017.
	\newblock ISSN 23752548.
	\newblock \doi{10.1126/sciadv.1701533}.
	\newblock URL \url{http://advances.sciencemag.org/}.
	
	\bibitem[Farazmand and Sapsis(2019{\natexlab{a}})]{farazmand2019b}
	M.~Farazmand and T.~P. Sapsis.
	\newblock Closed-loop adaptive control of extreme events in a turbulent flow.
	\newblock \emph{Phys. Rev. E}, 100:\penalty0 033110, 2019{\natexlab{a}}.
	\newblock \doi{10.1103/PhysRevE.100.033110}.
	
	\bibitem[Zhang et~al.(2020)Zhang, Xu, Liu, Kurths, and Grebogi]{Zhang2020}
	X.~Zhang, Y.~Xu, Q.~Liu, J.~Kurths, and C.~Grebogi.
	\newblock {Rate-dependent bifurcation dodging in a thermoacoustic system driven
		by colored noise (preprint)}.
	\newblock \emph{arXiv:2009.13126}, sep 2020.
	\newblock URL \url{https://arxiv.org/abs/2009.13126}.
	
	\bibitem[Farazmand and Sapsis(2019{\natexlab{b}})]{farazmand2019a}
	M.~Farazmand and T.~P. Sapsis.
	\newblock Extreme events: {M}echanisms and {P}rediction.
	\newblock \emph{Applied Mechanics Review}, 2019{\natexlab{b}}.
	\newblock \doi{10.1115/1.4042065}.
	
	\bibitem[Farazmand(2020)]{Farazmand2020}
	M.~Farazmand.
	\newblock {Mitigation of tipping point transitions by time-delay feedback
		control}.
	\newblock \emph{Chaos}, 30\penalty0 (013149), 2020.
	\newblock \doi{10.1063/1.5137825}.
	\newblock URL \url{https://doi.org/10.1063/1.5137825}.
	
	\bibitem[Horsthemke and Lefever(2006)]{Horsthemke2006}
	W.~Horsthemke and R.~Lefever.
	\newblock \emph{{Noise-Induced Transitions}}.
	\newblock Springer, 2nd edition, 2006.
	
	\bibitem[Pugachev and Sinitsyn(2001)]{Pugachev2001}
	V.S. Pugachev and I.N. Sinitsyn.
	\newblock \emph{{Stochastic Systems. Theory and Applications}}.
	\newblock World Scientific, 2001.
	
	\bibitem[H{\"{a}}nggi and Jung(1995)]{Hanggi1995}
	P.~H{\"{a}}nggi and P.~Jung.
	\newblock {Colored Noise in Dynamical Systems}.
	\newblock \emph{Advances in Chemical Physics}, 89:\penalty0 239--326, 1995.
	
	\bibitem[Mamis et~al.(2019)Mamis, Athanassoulis, and Kapelonis]{Mamis2019d}
	K.I. Mamis, G.~A. Athanassoulis, and Z.~G. Kapelonis.
	\newblock {A systematic path to non-Markovian dynamics: New response
		probability density function evolution equations under Gaussian coloured
		noise excitation}.
	\newblock \emph{Proc.R.Soc.A}, 471\penalty0 (20180837), 2019.
	
	\bibitem[{\O}ksendal(2003)]{Oksendal2003}
	B.~{\O}ksendal.
	\newblock \emph{{Stochastic Differential Equations. An Introduction with
			Applications}}.
	\newblock Springer, 6th edition, 2003.
	
	\bibitem[Risken(1996)]{Risken1996}
	H.~Risken.
	\newblock \emph{{The Fokker-Planck Equation. Methods of Solution and
			Applications}}.
	\newblock Springer-Verlag, 2nd edition, 1996.
	
	\bibitem[Gardiner(2004)]{Gardiner2004}
	C.W. Gardiner.
	\newblock \emph{{Handbook of Stochastic Methods for Physics, Chemistry and the
			Natural Sciences}}.
	\newblock Springer, 3rd edition, 2004.
	
	\bibitem[Roberts and Spanos(2003)]{Roberts2003}
	J.~B. Roberts and P.~D. Spanos.
	\newblock \emph{{Random Vibration and Statistical Linearization}}.
	\newblock Dover Publications, 2003.
	
	\bibitem[Spanos(1986)]{Spanos1986a}
	P.~T.D. Spanos.
	\newblock {Filter approaches to wave kinematics approximation}.
	\newblock \emph{Studies in Applied Mechanics}, 14\penalty0 (C):\penalty0
	459--473, jan 1986.
	\newblock ISSN 09225382.
	\newblock \doi{10.1016/B978-0-444-42665-9.50033-5}.
	
	\bibitem[Francescutto and Naito(2004)]{Francescutto2004}
	A.~Francescutto and S.~Naito.
	\newblock {Large amplitude rolling in a realistic sea}.
	\newblock \emph{International shipbuilding progress}, 51:\penalty0 221--235,
	2004.
	\newblock ISSN 0020868X.
	\newblock URL \url{http://iospress.metapress.com/index/D58TRBX995AFJRJL.pdf}.
	
	\bibitem[Scruggs et~al.(2013)Scruggs, Lattanzio, Taflanidis, and
	Cassidy]{Scruggs2013}
	J.~T. Scruggs, S.~M. Lattanzio, A.~A. Taflanidis, and I.~L. Cassidy.
	\newblock {Optimal causal control of a wave energy converter in a random sea}.
	\newblock \emph{Applied Ocean Research}, 42:\penalty0 1--15, aug 2013.
	\newblock ISSN 01411187.
	\newblock \doi{10.1016/j.apor.2013.03.004}.
	
	\bibitem[Chai et~al.(2015)Chai, Naess, and Leira]{Chai2015}
	W.~Chai, A.~Naess, and B.~J. Leira.
	\newblock {Filter models for prediction of stochastic ship roll response}.
	\newblock \emph{Probabilistic Engineering Mechanics}, 41:\penalty0 104--114,
	2015.
	\newblock ISSN 02668920.
	\newblock \doi{10.1016/j.probengmech.2015.06.002}.
	\newblock URL
	\url{http://linkinghub.elsevier.com/retrieve/pii/S0266892015300175}.
	
	\bibitem[Toral and Colet(2014)]{Toral2014}
	R.~Toral and P.~Colet.
	\newblock \emph{{Stochastic Numerical Methods. An Introduction for Students and
			Scientists}}.
	\newblock Wiley, 2014.
	
	\bibitem[Pyragas(1995)]{Pyragas1995}
	K.~Pyragas.
	\newblock {Control of chaos via extended delay feedback}.
	\newblock \emph{Physics Letters A}, 206\penalty0 (5-6):\penalty0 323--330, oct
	1995.
	\newblock ISSN 03759601.
	\newblock \doi{10.1016/0375-9601(95)00654-L}.
	
	\bibitem[Suresh and Chandrasekar(2018)]{Suresh2018}
	R.~Suresh and V.~K. Chandrasekar.
	\newblock {Influence of time-delay feedback on extreme events in a forced
		Li{\'{e}}nard system}.
	\newblock \emph{Physical Review E}, 98\penalty0 (5):\penalty0 052211, nov 2018.
	\newblock ISSN 24700053.
	\newblock \doi{10.1103/PhysRevE.98.052211}.
	\newblock URL
	\url{https://journals.aps.org/pre/abstract/10.1103/PhysRevE.98.052211}.
	
	\bibitem[Guillouzic et~al.(1999)Guillouzic, L'Heureux, and
	Longtin]{Longtin1999}
	S.~Guillouzic, I.~L'Heureux, and A.~Longtin.
	\newblock Small delay approximation of stochastic delay differential equations.
	\newblock \emph{Phys. Rev. E}, 59:\penalty0 3970--3982, Apr 1999.
	\newblock \doi{10.1103/PhysRevE.59.3970}.
	
	\bibitem[Sun(2006)]{Sun2006}
	J.Q. Sun.
	\newblock \emph{{Stochastic Dynamics and Control}}.
	\newblock Elsevier, 2006.
	
	\bibitem[Mamis and Athanassoulis(2016)]{Mamis2016e}
	K.I. Mamis and G.A. Athanassoulis.
	\newblock {Exact stationary solutions to Fokker-Planck-Kolmogorov equation for
		oscillators using a new splitting technique and a new class of stochastically
		equivalent systems}.
	\newblock \emph{Probabilistic Engineering Mechanics}, 45:\penalty0 22--30,
	2016.
	\newblock ISSN 18784275 02668920.
	\newblock \doi{10.1016/j.probengmech.2016.02.003}.
	
	\bibitem[Masud and Bergman(2005)]{Masud2005}
	A.~Masud and L.A. Bergman.
	\newblock {Solution of the Four Dimensional Fokker-Planck Equation: Still a
		Challenge}.
	\newblock In G.~Augusti, G.I. Schueller, and M.~Ciampoli, editors,
	\emph{ICOSSAR}, pages 1911--1916, Rotterdam, 2005.
	\newblock URL
	\url{http://scholar.google.com/scholar?hl=en{\&}btnG=Search{\&}q=intitle:Solution+of+the+Four+Dimensional+Fokker-Planck+Equation+:+Still+a+Challenge{\#}0}.
	
	\bibitem[Chen and Majda(2017)]{Chen2017}
	N.~Chen and A.~J. Majda.
	\newblock {Beating the curse of dimension with accurate statistics for the
		Fokker-Planck equation in complex turbulent systems}.
	\newblock \emph{Proceedings of the National Academy of Sciences of the United
		States of America}, 114\penalty0 (49):\penalty0 12864--12869, 2017.
	\newblock ISSN 10916490.
	\newblock \doi{10.1073/pnas.1717017114}.
	
	\bibitem[Xu et~al.(2020)Xu, Zhang, Li, Zhou, Liu, and Kurths]{Xu2020}
	Y.~Xu, H.~Zhang, Y.~Li, K.~Zhou, Q.~Liu, and J.~Kurths.
	\newblock {Solving Fokker-Planck equation using deep learning}.
	\newblock \emph{Chaos}, 30\penalty0 (1), 2020.
	\newblock ISSN 10541500.
	\newblock \doi{10.1063/1.5132840}.
	\newblock URL \url{https://doi.org/10.1063/1.5132840}.
	
	\bibitem[Fox(1987)]{Fox1987}
	R.~F. Fox.
	\newblock {Stochastic calculus in physics}.
	\newblock \emph{Journal of Statistical Physics}, 46\penalty0 (5-6):\penalty0
	1145--1157, 1987.
	\newblock ISSN 0022-4715.
	\newblock \doi{10.1007/BF01011160}.
	\newblock URL \url{http://link.springer.com/10.1007/BF01011160}.
	
	\bibitem[Faetti et~al.(1988)Faetti, Fronzoni, Grigolini, and
	Mannella]{Faetti1988}
	S.~Faetti, L.~Fronzoni, P.~Grigolini, and R.~Mannella.
	\newblock {The Projection Approach to the Fokker-Planck Equation. I. Colored
		Gaussian Noise}.
	\newblock \emph{Journal of Statistical Physics}, 52\penalty0 (3-4):\penalty0
	951--977, 1988.
	
	\bibitem[Peacock-L{\'{o}}pez et~al.(1988)Peacock-L{\'{o}}pez, West, and
	Lindenberg]{Peacock-Lopez1988}
	E.~Peacock-L{\'{o}}pez, B.~J. West, and K.~Lindenberg.
	\newblock {Relations among effective Fokker-Planck for systems driven by
		colored noise}.
	\newblock \emph{Physical Review A}, 37\penalty0 (9):\penalty0 3530--3535, 1988.
	
	\bibitem[H{\"{a}}nggi(1989)]{Hanggi1989}
	P.~H{\"{a}}nggi.
	\newblock {Colored noise in continuous dynamical systems: a functional calculus
		approach}.
	\newblock In F.~Moss and P.~V.E. McClintock, editors, \emph{Noise in nonlinear
		dynamical systems, vol. 1: Theory of continuous Fokker-Planck systems}, pages
	307--328. Cambridge University Press, 1989.
	
	\bibitem[Bianucci and Mannella(2020)]{Bianucci2020a}
	M.~Bianucci and R.~Mannella.
	\newblock {Optimal FPE for non-linear 1d-SDE. I: Additive Gaussian colored
		noise}.
	\newblock \emph{Journal of Physics Communications}, 4\penalty0 (10):\penalty0
	105019, nov 2020.
	\newblock ISSN 2399-6528.
	\newblock \doi{10.1088/2399-6528/abc54e}.
	\newblock URL
	\url{https://iopscience.iop.org/article/10.1088/2399-6528/abc54e}.
	
	\bibitem[H{\"{a}}nggi et~al.(1985)H{\"{a}}nggi, Mroczkowski, Moss, and
	McClintock]{Hanggi1985}
	P.~H{\"{a}}nggi, T.~J. Mroczkowski, F.~Moss, and P.~V.E. McClintock.
	\newblock {Bistability driven by colored noise: Theory and experiment}.
	\newblock \emph{Physical Review A}, 32\penalty0 (1):\penalty0 695--698, 1985.
	\newblock ISSN 10502947.
	\newblock \doi{10.1103/PhysRevA.32.695}.
	
	\bibitem[Frank(2005)]{Frank2005}
	T.~D. Frank.
	\newblock \emph{{Nonlinear Fokker-Planck Equations}}.
	\newblock Springer, 2005.
	
	\bibitem[{\"{O}}ttinger(1996)]{Ottinger1996}
	H.~C. {\"{O}}ttinger.
	\newblock \emph{{Stochastic Processes in Polymeric Fluids}}.
	\newblock Springer-Verlag, Berlin, Heidelberg, 1996.
	
	\bibitem[Sancho et~al.(1982)Sancho, {San Miguel}, Katz, and Gunton]{Sancho1982}
	J.M. Sancho, M.~{San Miguel}, S.L. Katz, and J.D. Gunton.
	\newblock {Analytical and numerical studies of multiplicative noise}.
	\newblock \emph{Physical Review A}, 26\penalty0 (3):\penalty0 1589--1609, 1982.
	
	\bibitem[Cao et~al.(2015)Cao, Zhang, and Karniadakis]{Cao2015}
	W.~Cao, Z.~Zhang, and G.~E. Karniadakis.
	\newblock {Numerical Methods for Stochastic Delay Differential Equations via
		the Wong-Zakai Approximation}.
	\newblock \emph{SIAM J. Sci. Comput.}, 37\penalty0 (1):\penalty0 295--318,
	2015.
	
	\bibitem[Wu et~al.(2009)Wu, Singh, and Xiao]{Wu2009}
	H.~Wu, S.~Singh, and M.~Xiao.
	\newblock {Multiplicative noise-induced probability distributions in
		three-level atomic optical bistability}.
	\newblock \emph{Physical Review A - Atomic, Molecular, and Optical Physics},
	79\penalty0 (2), 2009.
	\newblock ISSN 10502947.
	\newblock \doi{10.1103/PhysRevA.79.023835}.
	
	\bibitem[Geng et~al.(2020)Geng, Peters, Trichet, Malmir, Kolkowski, Smith, and
	Rodriguez]{Geng2020}
	Z.~Geng, K.~J.H. Peters, A.~A.P. Trichet, K.~Malmir, R.~Kolkowski, J.~M. Smith,
	and S.~R.K. Rodriguez.
	\newblock {Universal Scaling in the Dynamic Hysteresis, and Non-Markovian
		Dynamics, of a Tunable Optical Cavity}.
	\newblock \emph{Physical Review Letters}, 124\penalty0 (15):\penalty0 153603,
	2020.
	\newblock ISSN 10797114.
	\newblock \doi{10.1103/PhysRevLett.124.153603}.
	\newblock URL \url{https://doi.org/10.1103/PhysRevLett.124.153603}.
	
	\bibitem[Peters(2019)]{Peters2019}
	K.~J.~H. Peters.
	\newblock \emph{{Non-Markovian Stochastic Resonance in a Tunable Optical
			Microcavity}}.
	\newblock Master's thesis, Utrecht University, 2019.
	
	\bibitem[Mamis et~al.(2018)]{Mamis2018}
           K.I.~Mamis, G.A.~Athanassoulis, and K.E~Papadopoulos.
          \newblock{Generalized FPK equations corresponding to systems of nonlinear random differential    equations excited by colored noise. Revisitation and new directions}.
          \newblock \emph{Procedia Computer Science}, 136\penalty0 (C):\penalty0 164--173, 2018.
          \newblock \doi{10.1016/j.procs.2018.08.249}.
	
	\bibitem[Mamis(2020)]{Mamis2020}
	K.I. Mamis.
	\newblock \emph{{Probabilistic responses of dynamical systems subjected to
			Gaussian coloured noise excitation. Foundations of a non-Markovian theory}}.
	\newblock Phd thesis, National Technical University of Athens, 2020.
	\newblock URL
	\url{https://www.didaktorika.gr/eadd/handle/10442/47678?locale=en}.
	
	\bibitem[Soize(1994)]{Soize1994}
	C.~Soize.
	\newblock \emph{{The Fokker-Planck equation for stochastic dynamical systems
			and its explicit steady state solutions}}.
	\newblock World Scientific, 1994.
	
	\bibitem[Zhu et~al.(1992)Zhu, Cai, and Lin]{Zhu1992}
	W.~Q. Zhu, G.~Q. Cai, and Y.~K. Lin.
	\newblock {Stochastic Excited Hamiltonian Systems}.
	\newblock In \emph{Nonlinear Stochastic Mechanics IUTAM Symposium}, pages
	543--552, Turin, Italy, 1992.
	\newblock \doi{10.1007/978-3-642-84789-9_47}.
	
	\bibitem[Zhu and Huang(2001)]{Zhu2001}
	W.~Q. Zhu and Z.~L. Huang.
	\newblock {Exact stationary solutions of stochastically excited and dissipated
		partially integrable Hamiltonian systems}.
	\newblock \emph{International Journal of Non-Linear Mechanics}, 36:\penalty0
	39--48, 2001.
	
	\bibitem[Fox(1986)]{Fox1986a}
	R.~F. Fox.
	\newblock {Uniform convergence to an effective Fokker-Planck equation for
		weakly colored noise}.
	\newblock \emph{Physical Review A}, 34\penalty0 (5):\penalty0 4525, 1986.
	\newblock ISSN 10502947.
	\newblock \doi{10.1103/PhysRevA.34.4525}.
	\newblock URL \url{http://ukpmc.ac.uk/abstract/MED/9897829}.
	
\end{thebibliography}

\end{document}